\documentclass[11pt]{article}
\usepackage{amssymb}
\usepackage{amsmath}
\usepackage{epsfig}

\def\NZQ{\mathbb}               
\def\NN{{\NZQ N}}

\def\ZZ{{\NZQ Z}}
\def\RR{{\NZQ R}}

\def\N1{{\NZQ 1}}
\def\N0{{\NZQ 0}}
\def\bm#1{\mbox{$\boldsymbol{#1}$}}
\newcommand{\va}{\bm{a}}  \newcommand{\vb}{\bm{b}}
\newcommand{\vc}{\bm{c}}  
\newcommand{\ve}{\bm{e}}

  \newcommand{\vr}{\bm{r}}
  
\newcommand{\vu}{\bm{u}}  \newcommand{\vv}{\bm{v}}
\newcommand{\vw}{\bm{w}}  \newcommand{\vx}{\bm{x}}
\newcommand{\vy}{\bm{y}}

\newcommand{\vzero}{\mbox{\bf 0}} 

\def\opn#1#2{\def#1{\operatorname{#2}}} 
\opn\Ker{Ker}
\opn\Coker{Coker}
\opn\Im{Im}
\opn\Hom{Hom}
\opn\aff{aff}
\opn\conv{conv}
\opn\cone{cone}
\opn\relint{relint}
\opn\st{st}
\opn\lk{lk}
\opn\cn{cn}
\opn\core{core}
\opn\vol{vol}
\opn\link{link}
\opn\star{star}
\opn\supp{supp}
\opn\rank{rank}
\opn\diam{diam}
\def\pot#1#2{#1[\kern-0.28ex[#2]\kern-0.28ex]}

\opn\dirlim{\underrightarrow{\lim}}
\opn\inivlim{\underleftarrow{\lim}}
\def\Implies{\ifmmode\Longrightarrow \else
     \unskip\({}\Longrightarrow{}\)\ignorespaces\fi}
\def\implies{\ifmmode\Rightarrow \else
     \unskip\({}\Rightarrow{}\)\ignorespaces\fi}
\def\iff{\ifmmode\Longleftrightarrow \else
     \unskip\({}\Longleftrightarrow{}\)\ignorespaces\fi}

\let\:=\colon
\newtheorem{thm}{Theorem}[section]
\newtheorem{lem}[thm]{Lemma}

\newtheorem{prop}[thm]{Proposition}

\newtheorem{exmp}[thm]{Example}
\newtheorem{defn}[thm]{Definition}
\newtheorem{algo}[thm]{Algorithm}

\newtheorem{prob}[thm]{Problem}
\newenvironment{proof}{\medskip                    
\noindent{\it Proof:}}{\quad \(\square\)\medskip}  
\newenvironment{proof2}[1]{\medskip                    
\noindent{\it Proof of #1:}}{\hfill \(\square\)\medskip}  

\def\abs#1{\left| #1 \right| }
\def\any{\mbox{}^{\forall}}
\def\fitter{\mbox{}^{\exists}}

\newcommand{\arithdeg}[1]{arith\mbox{-}deg\left({#1}\right)}
\newcommand{\nidan}[2]{\genfrac{}{}{0pt}{}{{#1}}{{#2}}}
\newcommand{\lw}[1]{\smash{\lower2.ex\hbox{#1}}}
\def\T{^{\rm T}}

\begin{document}

\title{Dualistic computational algebraic analyses of primal and dual
  minimum cost flow problems on acyclic tournament graphs}
\author{
{\large \hspace*{-1ex} \(\mbox{Takayuki Ishizeki}^{\ast}\)}
\qquad
{\large \hspace*{-1ex} \(\mbox{Hiroki Nakayama}^{\dagger}\)}
\qquad
{\large \hspace*{-1ex} \(\mbox{Hiroshi Imai}^{\dagger}\)}
\\
{\tt \{ishizeki, nak-den, imai\}@is.s.u-tokyo.ac.jp}}

\maketitle

\footnotetext[1]{
Department of Information Science,
        Graduate School of Science, The University of Tokyo,
        7-3-1 Hongo, Bunkyo-ku, Tokyo 113-0033, Japan.}
\footnotetext[2]{
Department of Computer Science, Graduate School of Information Science
and Technology, The University of Tokyo, 7-3-1 Hongo, Bunkyo-ku, Tokyo 
113-0033, Japan.}

\vspace*{-1cm}
\begin{abstract}
  To integer programming problems, computational algebraic approaches
  using Gr{\"o}bner bases or standard pairs via the discreteness of
  toric ideals have been studied in recent years. Although these
  approaches have not given improved time complexity bound compared with
  existing methods for solving integer programming problems, these give 
  algebraic analysis of their structures. In this paper, we focus on the 
  case that the coefficient matrix is unimodular, especially on the
  primal and dual minimum cost flow problems, whose structure is rather
  well-known, but new structures can be revealed by our approach. We
  study the Gr{\"o}bner bases and standard pairs for unimodular
  programming, and give the maximum number of dual feasible bases in
  terms of the volume of polytopes. And for the minimum cost flow
  problems, we characterize reduced Gr{\"o}bner bases in terms of 
  graphs, and give bounds for the number of dual (resp. primal) feasible
  bases of the primal (resp. dual) problems: for the primal problems the
  minimum and the maximum are shown to be 1 and the Catalan number
  \(\frac{1}{d}\tbinom{2(d-1)}{d-1}\), while for the dual problems the
  lower bound is shown to be \(\Omega(2^{\lfloor d/6\rfloor})\). To
  analyze arithmetic degrees, we use two approaches: one is the relation
  between reduced Gr{\"o}bner bases and standard pairs, where the
  corresponding relation on the minimum cost flow --- between a subset
  of circuits and dual feasible bases --- has not been so clear, the
  other is the results in combinatorics related with toric ideals.
\end{abstract}

\section{Introduction}
Recently, some algebraic approaches to integer programming problems have 
been studied. The two main approaches are using {\em Gr{\"o}bner\/
  bases\/}~\cite{ConTra91} and {\em standard\/ pair\/
decompositions\/}~\cite{HosTho99StandardPair}. Although they neither
give improved complexity bounds compared with existing methods nor have
been demonstrated to solve hard practical instances which cannot be
handled by existing methods, these approaches themselves are very
interesting by applying computational algebraic methods to such hard
problems, and give algebraic analysis of structure of integer
programming
problems~\cite{ConTra91,HosTho99StandardPair,HosTho01,StuTho97,Tho95,Tho97,UrbWeiZie97,Zie99}.
For an ideal over a polynomial ring, the reduced Gr{\"o}bner basis and
the set of standard pairs are dual in a sense that the complement of the
monomials in the initial ideal, which is generated by initial terms of
the reduced Gr{\"o}bner basis, is the set of standard monomials, whose
nice decomposition is the standard pair decomposition. This kind of
duality may shed new light on duality in combinatorial optimization, and
by considering a nice subclass of integer programming problems where the
duality theorem holds, we might be able to obtain some complexity bounds
by making use of the characteristics of the subclass, which could not be
derived for general integer programming problems.

The problems whose coefficient matrices are unimodular form a nice subclass 
in a sense that the system \(yA\leq c\) becomes totally dual integral
(TDI). Then each standard pair corresponds to a dual feasible basis, and
the method using standard pairs is equivalent to calculate the reduced
cost for each basis\ (Theorem~\ref{thm:red_cost}). Thus, the number of
standard pairs, which is equal to that of dual feasible bases, gives the
complexity of this approach. Additionally, the maximum number of
standard pairs can be described by the normalized volume of another
matrix\ (Theorem~\ref{thm:MaxVolume}). 

Especially, the minimum cost flow problems form a well-known subclass of 
unimodular integer programming problems which can be solved in
polynomial time. A Gr{\"o}bner basis approach for the minimum cost flow
problems is a variant of cycle-canceling algorithm. In the case of the
strongly polynomial time
algorithms~\cite{GolTar89,ShiIwaMac00,SokAhuOrl00}, for any feasible
flow they choose polynomial size of negative-cost cycles (by the
selecting rules) from the set of negative-cost cycles in the residual
network, which may be of exponential size, as many as
possible. Similarly, the algorithm using Gr{\"o}bner basis calculates
the optimal flow by augmenting flows along the negative-cost cycles
which correspond to the elements of Gr{\"o}bner basis. Thus the
cardinalities of reduced Gr{\"o}bner bases may give some time bound for
this algorithm. On the other hand, a standard pair approach for the
minimum cost flow problems first finds the set of standard pairs, and
solves linear system of equations for each standard pair until an
integer and non-negative solution is obtained.
For a network optimization problem, the duality between the reduced
Gr{\"o}bner basis and the set of standard pairs corresponds to the
relation between circuits and dual feasible co-trees, dually, cutsets
and primal feasible trees. Since such a relation has not been so clear,
the computational algebraic duality may be interesting method for the
analysis of network problems.

This paper is organized as follows. In Section 2, reduced Gr{\"o}bner
bases and standard pairs are defined, and their relations with integer
programming problems, regular triangulations and dual polyhedra are
introduced. The case that the coefficient matrix is unimodular is
studied in Section 3. The maximum arithmetic degree (i.e. the maximum
number of dual feasible bases) is shown to be the normalized volume of a
polytope defined by homogenizing the coefficient matrix\ 
(Theorem~\ref{thm:MaxVolume}). In Section 4, we study the Gr{\"o}bner
bases and standard pairs on the primal minimum cost flow problems on
acyclic tournament graphs with \(d\) vertices. We show that three types
of reduced Gr{\"o}bner bases can be characterized in terms of the
circuits\ (Theorem~\ref{thm:grob_p_type1},\ \ref{thm:grob_p_type2},\
\ref{thm:grob_p_type3}). These examples give the minimum and the maximum
number of dual feasible bases of the minimum cost flow problems: the
minimum is 1\ (Theorem~\ref{thm:min_arith_degree}) and the maximum is
\(\frac{1}{d}\tbinom{2(d-1)}{d-1}\)\
(Theorem~\ref{thm:max_arith_degree}). This maximum is shown using the
result in Section 3 and the result about the hypergeometric system on
unipotent matrices and related polytope~\cite{GelGraPos96}. In Section
5, we study the dual minimum cost flow problems. One reduced Gr{\"o}bner
basis is characterized in terms of cutsets\
(Theorem~\ref{thm:grob_d}). We also show that the lower bound for the
number of primal feasible bases of the minimum cost flow problems is
\(\Omega(2^{\lfloor d/6\rfloor})\). 

\begin{table}[htbp]
  \begin{center}
    \begin{tabular}{|c|c@{\quad\vrule width1.0pt\quad}l|l|}
      \hline
      \multicolumn{2}{|c@{\quad\vrule width1.0pt\quad}}{ } & Gr{\"o}bner
      basis~\cite{ConTra91} & Standard
      pair~\cite{HosTho99StandardPair}\\
      \noalign{\hrule height 1.0pt}
      & Term on graph & Set of circuits & Set of spanning trees\\
      \cline{2-4}
      & \lw{Algorithm} & Variant of & Enumeration of\\
      Primal & & cycle-canceling & dual feasible bases\\
      \cline{2-4}
      & On acyclic tournament& \(\min:\ d(d-1)/2\) & \(\min:\ 1\)\\
      & graph with \(d\) vertices & \(\max:\ \) ? & \(\max:\
      \frac{1}{d}\tbinom{2(d-1)}{d-1}\)\\
      \hline
      & Term on graph & Set of cutsets & Set of co-trees\\
      \cline{2-4}
      & \lw{Algorithm} & Variant of & Enumeration of\\
      Dual & & cutset-canceling & primal feasible bases\\
      \cline{2-4}
      & On acyclic tournament & \(\min:\ d-1\) & \lw{Lower bound
        \(\Omega(2^{\lfloor d/6\rfloor})\)}\\
      & graph with \(d\) vertices & \(\max:\ \) ? & \\
      \hline
    \end{tabular}
  \end{center}
  \caption{Dual algebraic approaches for primal and dual minimum cost
    flow problems}
  \label{fig:dual_approaches}
\end{table}

\section{Toric ideals and Gr{\"o}bner bases}
For a matrix \(A\in\ZZ^{d\times n}\) and a cost vector \(\vc\in\RR^n\),
let \(IP_{A,\vc}\) be the family of integer programming problems
\(IP_{A,\vc}(\vb):=minimize\ \{\vc\cdot\vx | A\vx=\vb,\ \vx\in\NN^n\}\)
as \(\vb\) varies in \(\{A\vu | \vu\in\NN^n\}\subseteq\ZZ^d\) (\(\NN\)
is the set of non-negative integers). The cost vector \(\vc\) is called
{\em generic\/} if each program in \(IP_{A,\vc}\) has the unique optimal
solution.

Let \(k\) be a field and \(k[\vx]:=k[x_1,\ldots,x_n]\) the polynomial
ring. For an exponent vector \(\va=(a_1,\ldots,a_n)\in\NN^n\), we denote 
\(\vx^{\va}:=x_1^{a_1}x_2^{a_2}\cdots x_n^{a_n}\). A total order on
monomials in \(k[\vx]\) is a {\em term\/ order\/} if 1 is the unique
minimal element, and \(\vx^{\vu}\succ\vx^{\vv}\) implies
\(\vx^{\vu+\vw}\succ\vx^{\vv+\vw}\) for all \(\vu,\vv,\vw\in\NN^n\). For
a fixed term order \(\succ\), the {\em refinement\/} \(\succ_{\vc}\) of
\(\vc\) by \(\succ\) is a total order such that
\(\vx^{\vu}\succ_{\vc}\vx^{\vv}\) if either \(\vc\cdot\vu>\vc\cdot\vv\) or
``\(\vc\cdot\vu=\vc\cdot\vv\) and \(\vx^{\vu}\succ\vx^{\vv}\)''
holds. If \(\vc\geq \vzero\), then \(\succ_{\vc}\) becomes a term
order.

The {\em toric\/ ideal\/} \(I_A\) of \(A\) is a binomial ideal
\(I_A:=\langle\vx^{\vu}-\vx^{\vv}\ |\  A\vu=A\vv,\
\vu,\vv\in\NN^n\rangle\). For any \(f\in I_A\), the {\em initial\/
  term\/} \(in_{\succ_{\vc}}(f)\) of \(f\) is the largest term in \(f\)
with respect to \(\succ_{\vc}\). Then we define the {\em initial\/
  ideal\/} \(in_{\succ_{\vc}}(I_A)\) of \(I_A\) as
\(in_{\succ_{\vc}}(I_A):=\langle in_{\succ_{\vc}}(f)\ |\  f\in
I_A\rangle\).

\subsection{Gr{\"o}bner bases and Conti-Traverso algorithm}
A finite subset \({\mathcal G}_{\succ_{\vc}}=\{g_1,\ldots,g_s\}\subseteq
I_A\) is a {\em Gr{\"o}bner\/ basis\/} for \(I_A\) with respect to
\(\succ_{\vc}\) if \(in_{\succ_{\vc}}(I_A)=\langle
in_{\succ_{\vc}}(g_1),\ldots, in_{\succ_{\vc}}(g_s)\rangle\).  In
addition, Gr{\"o}bner basis \({\mathcal G}_{\succ_{\vc}}\) is {\em
  reduced\/} if \({\mathcal G}_{\succ_{\vc}}\) satisfies that (i) for
any \(i\), the coefficient of \(in_{\succ_{\vc}}(g_i)\) is \(1\), and
(ii) for any \(i\), any term of \(g_i\) is not divisible by
\(in_{\succ_{\vc}}(g_j)\ (i\neq j)\). If \(\succ_{\vc}\) is a term
order, then the reduced Gr{\"o}bner basis \({\mathcal G}_{\succ_{\vc}}\) 
exists uniquely, and is calculated by Buchberger algorithm (see
\cite{CoxLitOsh96}). Any Gr{\"o}bner basis for \(I_A\) is a basis of
\(I_A\)~\cite{CoxLitOsh96}.

\(I_A\) is called {\em homogeneous\/} with respect to the positive
grading \(\deg(x_i)=d_i>0\ (i=1,\ldots,n)\) if, for any
\(f=f_1+f_2+\cdots+f_m\in I_A\) (\(f_i\) is the homogeneous component of
degree \(i\) in \(f\)), \(f_i\in I_A\) for any \(i\). Then \(I_A\) is
homogeneous if and only if \(I_A\) is generated by homogeneous
polynomials~\cite{CoxLitOsh96}.

\begin{prop}[\cite{Stu95}]
  If \(I_A\) is a homogeneous with respect to some positive grading
  \(\deg(x_i)=d_i>0\), then \(\succ_{\vc}\) becomes a term order for any 
  \(\vc\in\RR^n\setminus\{\vzero\}\), and the reduced Gr{\"o}bner basis
  \({\mathcal G}_{\succ_{\vc}}\) exists.
  \label{prop:grad}
\end{prop}

The {\em support\/} \(\supp(\vu)\) of a vector \(\vu\) is the index set
\(\{i\ |\  u_i\neq 0\}\). Any \(\vu\in\ZZ^n\) can be written uniquely as
\(\vu=\vu^{+}-\vu^{-}\) where \(\vu^{+}, \vu^{-}\in\NN^n\) and have
disjoint support. Then \({\mathcal G}_{\succ_{\vc}}\) can be written as
\({\mathcal G}_{\succ_{\vc}}=\{\vx^{\vu_1^{+}}-\vx^{\vu_1^{-}},\ldots,
\vx^{\vu_p^{+}}-\vx^{\vu_p^{-}}\}\) for some finite
\(\vu_1,\ldots,\vu_p\in\ker(A)\cap\ZZ^n\)~\cite{Stu95}.

\begin{exmp}
  Let \(\displaystyle{A=\left(\begin{array}{ccc} 1 & 1 & 0\\ -1 & 0 &
        1\end{array}\right)}\) and consider the minimum cost flow
  problem \(IP_{A,\vc}(\vb)=minimize\{\vc\cdot\vx\ |\  A\vx=\vb,
  \vx=(x_{1,2},x_{1,3},x_{2,3})\in\NN^3\}.\)
  \begin{figure}[h]
    \begin{center}
      \epsfig{file=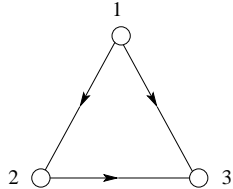,width=3cm}
      \caption{Acyclic tournament graph with 3 vertices.}
      \label{fig:K3}
    \end{center}
  \end{figure}
  Then the toric ideal is \(I_A=\langle
  x_{1,2}x_{2,3}-x_{1,3}\rangle\). If
  \(\vc=(c_{1,2},c_{1,3},c_{2,3})=(3,1,2)\), then
  \(in_{\vc}(I_A)=\langle x_{1,2}x_{2,3}\rangle\) and reduced
  Gr{\"o}bner basis is \({\mathcal
    G}_{\succ_{\vc}}=\{x_{1,2}x_{2,3}-x_{1,3}\}\).
  \label{exmp:K3}
\end{exmp}

In the rest of this paper, we consider a cost vector \(\vc\) which
\(\succ_{\vc}\) becomes a term order for some term order \(\succ\). Let
\(IP_{A,\succ_{\vc}}(\vb)\) be the problem to find the unique minimal
element in \(\{\vx\in\NN^n\ |\  A\vx=\vb\}\) with respect to
\(\succ_{\vc}\). Then the solution \(\vu\) of
\(IP_{A,\succ_{\vc}}(\vb)\) is one of the optimal solutions of
\(IP_{A,\vc}(\vb)\). Conti and Traverso~\cite{ConTra91} introduced an
algorithm based on a Gr{\"o}bner basis to solve
\(IP_{A,\succ_{\vc}}(\vb)\). We describe the condensed version of the
Conti-Traverso Algorithm (see~\cite{Stu95}), which is useful for
highlighting the main computational step involved. The {\em normal\/
  form\/} of \(f\in k[\vx]\) by the reduced Gr{\"o}bner basis
\({\mathcal G}\) is the unique remainder obtained upon dividing \(f\) by 
\({\mathcal G}\). 

\begin{algo}[Conti-Traverso Algorithm~\cite{ConTra91,Stu95}]
 \begin{tabbing}
  \\
  \quad\={\bf 1.}\ \= Compute the reduced Gr{\"o}bner basis \({\mathcal
    G}_{\succ_{\vc}}\) of \(I_A\) with respect to \(\succ_{\vc}\).\\
  \>{\bf 2.}\ \= For any feasible solution \(\vv\) of
  \(IP_{A,\vc}(\vb)\), compute the normal form \(\vx^{\vu}\) of
  \(\vx^{\vv}\) by \({\mathcal G}_{\succ_{\vc}}\).\\
  \>{\bf 3.}\ Output \(\vu\). \(\vu\) is the optimal solution of
  \(IP_{A,\succ_{\vc}}(\vb)\).
 \end{tabbing}
 \label{alg:CT1}
\end{algo}

Thus reduced Gr{\"o}bner basis \({\mathcal G}_{\succ_{\vc}}\) is a test
set for \(IP_{A,\succ_{\vc}}\)~\cite{Tho95,Tho97}. 

\setcounter{thm}{1}
\begin{exmp}[continued]
  Let \(\vb=(4,5,-9)\). For a feasible solution \((4,0,9)\), the normal
  form of \(x_{1,2}^4x_{2,3}^9\) by \({\mathcal G}_{\succ_{\vc}}\) is
  \(x_{1,3}^4x_{2,3}^5\). Thus the optimal solution of
  \(IP_{A,\succ_{\vc}}(\vb)\) is \((0,4,5)\).
\end{exmp}

\setcounter{thm}{3}
\subsection{Standard pair decompositions}
Let \([n]:=\{1,\ldots,n\}\). For a monomial \(\vx^{\va}\in k[\vx]\) and
an index set \(\sigma\subseteq [n]\), \((\vx^{\va},\sigma)\) is a {\em
  standard\/ pair\/} of \(in_{\succ_{\vc}}(I_A)\) if (i)
\(\supp(\va)\cap\sigma=\emptyset\), (ii) every monomial in
\(\vx^{\va}\cdot k[x_j\ |\  j\in\sigma]:=\{\vx^{\va}\cdot f\ |\ f\in
k[x_j\ |\ j\in\sigma]\}\) is not an element of
\(in_{\succ_{\vc}}(I_A)\), and (iii) there exists no other
\((\vx^{\va'},\sigma')\), which satisfies {\rm (i)} and {\rm (ii)}, such 
that \(\vx^{\va'}\) divides \(\vx^{\va}\) and
\(\supp(\vx^{\va}/\vx^{\va'})\cup\sigma\subseteq\sigma'\). We denote
\(S(in_{\succ_{\vc}}(I_A))\) the set of all standard pairs of
\(in_{\succ_{\vc}}(I_A)\). We use the same \((\vx^{\va},\sigma)\) to
denote the set of all monomials in \(\vx^{\va}\cdot k[x_j\ |\
j\in\sigma]\). Then the above condition (iii) says that
\((\vx^{\va},\sigma)\not\subset (\vx^{\va'},\sigma')\) for any other
\((\vx^{\va'},\sigma')\) which satisfies the condition (i) and (ii).
The standard pairs of \(in_{\succ_{\vc}}(I_A)\) induce a
unique covering for the set of standard monomials of
\(in_{\succ_{\vc}}(I_A)\), which we call the {\em standard\/ pair\/
  decomposition\/} of
\(in_{\succ_{\vc}}(I_A)\). \(\abs{S(in_{\succ_{\vc}}(I_A))}\) is called
the {\em arithmetic\/ degree\/} of \(in_{\succ_{\vc}}(I_A)\) and denoted 
by \(\arithdeg{in_{\succ_{\vc}}(I_A)}\)~\cite{StuTruVog95}.

\setcounter{thm}{1}
\begin{exmp}[continued]
  For \(\vc=(3,1,2)\),  the standard pairs of \(in_{(3,1,2)}(I_A)\) are
  \((1,\{(1,2),(1,3)\})\) and \((1,\{(1,3),(2,3)\})\), thus the
  arithmetic degree of \(in_{(3,1,2)}(I_A)\) is {\rm 2}. On the other
  hand, for \(\vc=(1,4,2)\), the standard pair of \(in_{(1,4,2)}(I_A)\)
  is \((1,\{(1,2),(2,3)\})\), thus the arithmetic degree of
  \(in_{(1,4,2)}(I_A)\) is {\rm 1}.
  \begin{figure}[h]
    \begin{center}
      \epsfig{file=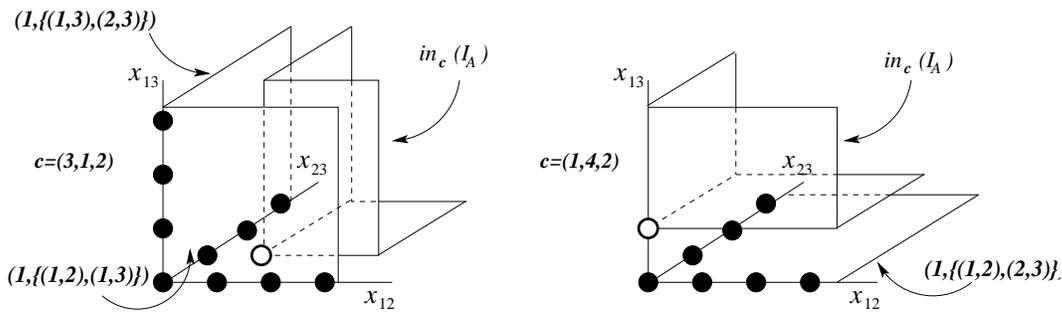,width=14cm}
      \caption{Two types of standard pair decompositions. A point
        \((p,q,r)\) in the figure corresponds to the monomial
        \(x_{1,2}^px_{1,3}^qx_{2,3}^r\).}
      \label{fig:std_pair_K3}
    \end{center}
  \end{figure}
\end{exmp}

Let \(\vc\) be a generic cost vector. Then
\(in_{\succ_{\vc}}(I_A)=in_{\vc}(I_A)\). Let \(\{\va_1,\ldots,\va_n\}\)
be the column vectors of \(A\) and \(\cone(A)\) the cone generated by
\(\va_1,\ldots,\va_n\). For \(\sigma\subseteq [n]\), we denote
\(A_{\sigma}\) for the submatrix of \(A\) whose columns are indexed by
\(\sigma\). For a cost vector \(\vc\), we define the {\em regular\/
  triangulation\/} \(\Delta_{\vc}\) of \(\cone(A)\) as follows:
\(\cone(A_{\sigma})\) is a face of \(\Delta_{\vc}\) if and only if there 
exists a vector \(\vy\in\RR^d\) such that \(\vy\cdot\va_j=c_j\
(j\in\sigma)\) and \(\vy\cdot\va_j<c_j\ (j\notin\sigma)\). If
\(\cone(A_{\sigma})\) is a face of \(\Delta_{\vc}\), \(\sigma\) also is
called a {\em face\/} of \(\Delta_{\vc}\). The genericity of \(\vc\)
implies that \(\Delta_{\vc}\) is in fact a triangulation (i.e. each face 
of \(\Delta_{\vc}\) is simplicial)~\cite{StuTho97}.

\setcounter{thm}{3}
\begin{lem}[\cite{Stu95,StuTruVog95}]
  \mbox{}
  \begin{enumerate}
    \setlength{\itemsep}{0pt}
  \item If \(in_{\vc}(I_A)\) has \((*,\sigma)\) as a standard pair, then
    \(\sigma\) is a face of \(\Delta_{\vc}\).
  \item \(in_{\vc}(I_A)\) has \((1,\sigma)\) as a standard pair if and
    only if \(\sigma\) is a maximal face of \(\Delta_{\vc}\).
  \item If \(\va_1,\ldots,\va_n\) span an affine hyperplane, then
    \(\Delta_{\vc}\) is the same as the regular triangulation of
    \(\conv(A)\) with respect to \(\vc\), and the  number of standard
    pairs \((*,\sigma)\) for a maximal face \(\sigma\) of
    \(\Delta_{\vc}\) equals the normalized volume of \(\sigma\) in
    \(\Delta_{\vc}\).
  \end{enumerate}
  \label{lem:Stu_and_STV}
\end{lem}

When vertices of \(\conv(A)\) are in the \(m\)-dimensional lattice
\(L\simeq\ZZ^m\), we define the {\em normalized\/ volume\/} of a maximal 
face \(\sigma\) of \(\Delta_{\vc}\) by the volume of \(\sigma\) with the 
normalization that the volume of the convex hull of
\(\vzero,\ve_1,\ldots,\ve_m\) is 1. Here, \(\{\ve_i\}_{1\leq i\leq m}\) are
the basis of the lattice \(L\).

For a polyhedron \(P\subset\RR^n\) and a face \(F\) of \(P\), the {\em
  normal\/ cone\/} of \(F\) at \(P\) is the cone
\(N_{P}(F):=\{\omega\in\RR^n\ |\ \omega\cdot x'\geq\omega\cdot x\
\mbox{for all}\ x'\in F\ \mbox{and}\ x\in P\}\). The set of normal cones
for all faces of \(P\) is called the {\em normal\/ fan\/} of \(P\).

\begin{lem}[\cite{HosTho01}]
  \(\Delta_{\vc}\) is the normal fan of the polyhedron
  \(P_{\vc}:=\{\vy\in\RR^d\ |\ \vy A\leq\vc\}\).
  \label{lem:Tri_and_Dual}
\end{lem}

We remark that \(P_{\vc}\) is the dual polyhedron for the linear
relaxation problem \(LP_{A,\vc}(\vb) := minimize\)\\ \(\left\{\vc\cdot\vx\ |\
  A\vx=\vb,\vx\geq\vzero\right\}\)
of \(IP_{A,\vc}(\vb)\). When \(A\) is row-full rank, this lemma shows
that there is one-to-one correspondence between the dual feasible bases
of \(LP_{A,\vc}(\vb)\) and the maximal faces of \(\Delta_{\vc}\).

\setcounter{thm}{1}
\begin{exmp}[continued]
  For \(\vc=(3,1,2)\),
  \(\Delta_{(3,1,2)}=\{\{1,2\},\{2,3\},\{1\},\{2\},\{3\},\emptyset\}\).
  \begin{figure}[h]
    \begin{center}
      \epsfig{file=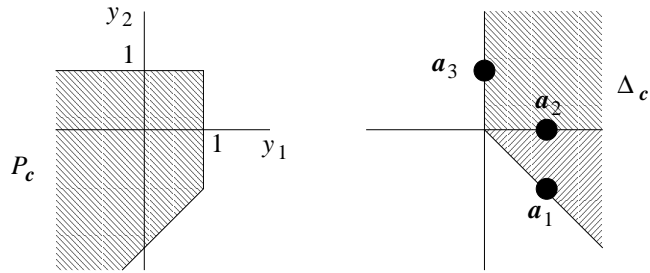,width=8.5cm}
      \caption{Dual polyhedron \(P_{(3,1,2)}\) and regular triangulation
        \(\Delta_{(3,1,2)}\).}
      \label{fig:dual_polyhedron_K3}
    \end{center}
  \end{figure}
\end{exmp}

Let \(\vu\) be the optimal solution to \(IP_{A,\vc}(\vb)\). Since
standard pairs cover \(in_{\vc}(I_A)\), \(\vx^{\vu}\) is covered by some
standard pair \((\vx^{\va},\sigma)\). Thus
\(\vu=\va+\sum_{i\in\sigma}k_i\ve_i\) for some non-negative integers
\(\{k_i\}_{i\in\sigma}\), and
\(\vb=A\vu=A\left(\va+\sum_{i\in\sigma}k_i\ve_i\right)=A\va+\sum_{i\in\sigma}k_i\va_i\).
Lemma~\ref{lem:Stu_and_STV} implies that \(\{\va_i\}_{i\in\sigma}\) are
linearly independent. Therefore \(\{k_i\}_{i\in\sigma}\) is the unique
solution to the linear system \(\sum_{i\in\sigma}x_i\va_i=\vb-A\va\). 
This observation induces an algorithm to solve \(IP_{A,\vc}(\vb)\) using
the standard pair decomposition of \(in_{\vc}(I_A)\).

\setcounter{thm}{5}
\begin{algo}
  {\bf (Solving \(IP_{A,\vc}(\vb)\) using \(S(in_{\vc}(I_A))\))}
 \begin{description}
  \setlength{\itemsep}{0pt}
  \item[(i)] For \((\vx^{\va},\sigma)\in S(in_{\vc}(I_A))\), solve the
    linear system \(\sum_{i\in\sigma}x_i\va_i=\vb-A\va\). Let
    \(\{k_i\}_{i\in\sigma}\) be the solution.
  \item[(ii)] If \(\{k_i\}_{i\in\sigma}\) are both integral and
    non-negative, output \(\va+\sum_{i\in\sigma}k_i\ve_i\) as the
    optimal solution. Otherwise, repeat {\rm (i)} for another standard
    pair.
 \end{description}
 \label{algo:IP_Standard_pair}
\end{algo}

This algorithm solves at most \(\arithdeg{in_{\vc}(I_A)}\)-many linear
systems of equations. Therefore \(\arithdeg{in_{\vc}(I_A)}\) is a
measure of the complexity of \(IP_{A,\vc}\).

\section{Standard pairs for unimodular programming}
\label{subsec:dual}
Let \(A\in\ZZ^{d\times n}\) be row-full rank and unimodular, i.e. each
non-zero maximal minor is \(\pm k\) for some \(k\in\NN\). Then
\(in_{\vc}(I_A)\) is minimally generated by square-free monomials for
any \(\vc\)~\cite{Stu95}, and all standard pairs are obtained from all
maximal faces of \(\Delta_{\vc}\). 

\begin{lem}[\cite{HosTho01}]
  Let \(\{m_1,\ldots,m_s\}\) be the minimal generators of
  \(in_{\vc}(I_A)\). If \(m_1,\ldots,m_s\) are all square-free then
  \(S(in_{\vc}(I_A))=\{(1,\sigma)\ |\ \sigma\ \mbox{is the maximal faces
    of}\ \Delta_{\vc}\}\). 
  \label{lem:square_free_std_pairs}
\end{lem}

For a matrix \(A\in\ZZ^{d\times n}\), the {\em homogenized\/ matrix\/}
\(A'\in\ZZ^{(d+1)\times(n+1)}\) of \(A\) is
\begin{equation}
  A':=\left(\begin{array}{cc}
    \begin{array}{cccc}
      1 & 1 & \cdots & 1
    \end{array} & 1\\
    A & \vzero
    \end{array}\right) = \left(\begin{array}{ccccc}
      1 & 1 & \cdots & 1 & 1\\
      \va_1 & \va_2 & \cdots & \va_n & \vzero
    \end{array}\right). \label{eq:homog_matrix}
\end{equation}
Let \(\va'_i=\tbinom{1}{\va_i}\) for \(1\leq i\leq n\) and \(\va'_{n+1}\) be
the \((n+1)\)-th column vector of \(A'\). We remark that
\(\va'_1,\ldots,\va'_n,\va'_{n+1}\) span an affine hyperplane.

We define another family \(IP_{A',(\vc,0)}\) of integer programming
problem
\[
 IP_{A',(\vc,0)}(\vb,\beta):=minimize\ \left\{\vc\cdot\vx\ \left|\
    A'\tbinom{\vx}{x_{n+1}}=\tbinom{\beta}{\vb},
    \tbinom{\vx}{x_{n+1}}\in\NN^{n+1}\right. \right\}
\]
as \(\tbinom{\beta}{\vb}\) varies in \(\{A'\vu\ |\
\vu\in\NN^{n+1}\}\). We remark that \((\vc,0)\) is generic if \(\vc\) is 
generic. 

The next proposition is due to Sturmfels et al.~\cite{StuTruVog95} for
general ideals. We give another proof for the case of toric ideal.

\begin{prop}[\cite{StuTruVog95}]
  \((\vx^{\va},\sigma)\in S(in_{\vc}(I_A))\ (\vx^{\va}\in
  k[\vx],\ \sigma\subseteq [n])\) if and only if \(\left(\vx^{\va},
    \sigma\cup\{n+1\}\right)\in S(in_{(\vc,0)}(I_{A'}))\).
  \label{prop:homog_std_pair}
\end{prop}

\begin{proof}
  We first show that any monomial in \((\vx^{\va},\sigma)\) is standard
  for \(in_{\vc}(I_A)\) if and only if any monomial in
  \((\vx^{\va},\sigma\cup\{n+1\})\) is standard for
  \(in_{(\vc,0)}(I_{A'})\). Suppose that any monomial in
  \((\vx^{\va},\sigma)\) is standard for \(in_{\vc}(I_A)\) and choose any
  \(\vx^{\vu}x_{n+1}^k\in\left(\vx^{\va},\sigma\cup\{n+1\}\right)\).
  If there exist any other \(\tbinom{\vv}{l}\in\NN^{n+1}\) such that
  \(A'\tbinom{\vv}{l}=A'\tbinom{\vu}{k}\) and
  \(\tbinom{\vv}{l}\neq\tbinom{\vu}{k}\), then \(A\vu=A\vv\), and
  \(\vc\cdot\vu<\vc\cdot\vv\) since \(\vx^{\vu}\notin
  in_{\vc}(I_A)\). Therefore, \(\tbinom{\vu}{k}\) is the optimal
  solution to \(IP_{A',(\vc,0)}(A\vu,\sum_{i=1}^{n}u_i+k)\). If there
  does not exist such \(\tbinom{\vv}{l}\), then clearly
  \(\tbinom{\vu}{k}\) is the optimal for
  \(IP_{A',(\vc,0)}(A\vu,\sum_{i=1}^{n}u_i+k)\). This shows that any
  monomial in \((\vx^{\va},\sigma\cup\{n+1\})\) is standard for
  \(in_{(\vc,0)}(I_{A'})\).
  
  Conversely, suppose that any monomial in
  \((\vx^{\va},\sigma\cup\{n+1\})\) is standard for \(in_{(\vc,0)}(I_{A'})\) 
  and choose any \(\vx^{\vu}\in(\vx^{\va},\sigma)\). If there exists some
  \(\vv\in\NN^n\) such that \(A\vv=A\vu\), then
  \(A'\tbinom{\vu}{p}=A'\tbinom{\vv}{q}\) for any non-negative integers
  \(p,\ q\) such that \(p-q=\sum_{i=1}^{n}v_i-\sum_{i=1}^{n}u_i\). Since
  \(\vx^{\vu}x_{n+1}^p\in\left(\vx^{\va},\sigma\cup\{n+1\}\right)\), 
  \((\vc,0)\cdot\tbinom{\vu}{p}<(\vc,0)\cdot\tbinom{\vv}{q}\), which implies
  that \(\vc\cdot\vu<\vc\cdot\vv\). Therefore, \(\vu\) is the optimal
  solution to \(IP_{A,\vc}(A\vu)\). If there does not exist such \(\vv\),
  then clearly \(\vu\) is the optimal for \(IP_{A,\vc}(A\vu)\). Thus any
  monomial in \((\vx^{\va},\sigma)\) is standard for \(in_{\vc}(I_A)\).

  Let \(\left(\vx^{\va},\sigma\cup\{n+1\}\right)\in
  S(in_{(\vc,0)}(I_{A'}))\). If
  \((\vx^{\va},\sigma)\subset(\vx^{\va'},\tau)\) for any other
  \((\vx^{\va'},\tau)\) which satisfies {\rm (i)} and {\rm (ii)} in the
  definition of standard pairs for \(in_{\vc}(I_A)\), then
  \(\left(\vx^{\va},\sigma\cup\{n+1\}\right)\) must be contained in
  \(\left(\vx^{\va'},\tau\cup\{n+1\}\right)\), which contradicts the
  assumption. Thus \((\vx^{\va},\sigma)\in S(in_{\vc}(I_A))\). On the
  other hand, if \(\left(\vx^{\va},\sigma\cup\{n+1\}\right)\notin
  S(in_{(\vc,0)}(I_{A'}))\), then there exists some
  \(\left(\vx^{\va'}x_{n+1}^k,\tau'\right)\) which contains
  \(\left(\vx^{\va},\sigma\cup\{n+1\}\right)\) and
  \(\left(\vx^{\va'}x_{n+1}^k,\tau'\right)\) satisfies {\rm (i)} and
  {\rm (ii)} in the definition of standard pairs for
  \(in_{(\vc,0)}(I_{A'})\). Then \(n+1\in\tau'\), and therefore
  \(k=0\). Therefore, \((\vx^{\va'},\tau'\setminus\{n+1\})\) contains
  \((\vx^{\va},\sigma)\) and satisfies {\rm (i)} and {\rm (ii)} in the
  definition of standard pairs. Thus \((\vx^{\va},\sigma)\notin
  S(in_{\vc}(I_A))\). This completes the proof.
\end{proof}

\setcounter{section}{2}
\setcounter{thm}{1}
\begin{exmp}[continued]
  For this \(A\), enlarged matrix \(A'\) is
  \[
   A'=\left(\begin{array}{cccc}
      1 & 1 & 1 & 1\\
      1 & 1 & 0 & 0\\
      -1 & 0 & 1 & 0
    \end{array}\right).\]
  We consider \(I_{A'}\subset k[x_{1,2},x_{1,3},x_{2,3},x_4]\). For
  \(\vc=(3,1,2)\), the standard pairs of \(in_{(3,1,2,0)}(I_{A'})\) are
  \((1,\{(1,2),(1,3),4\})\) and \((1,\{(1,3),(2,3),4\})\). On the other
  hand, for \(\vc=(1,4,2)\), the standard pairs of
  \(in_{(1,4,2,0)}(I_{A'})\) are \((1,\{(1,2),(1,3),(2,3)\})\) and
  \((1,\{(1,2),(2,3),4\})\). In this case, the only standard pair
  \((1,\{(1,2),(2,3),4\})\) satisfies the condition in
  Proposition~\ref{prop:homog_std_pair}, which corresponds to the
  standard pair \((1,\{(1,2),(2,3)\})\) of \(in_{(1,4,2)}(I_A)\).
\end{exmp}

Since \(\va'_1,\ldots,\va'_{n+1}\) span an affine hyperplane, the
normalized volume of \(\conv(A')\) gives the number of standard pairs of
\(in_{(\vc,k)}(I_{A'})\) which correspond to maximal faces of
\(\Delta'_{(\vc,k)}\) by Lemma~\ref{lem:Stu_and_STV}\ {\rm (iii)}.

\setcounter{section}{3}
\setcounter{thm}{2}
\begin{thm}[\cite{HibOhs02}]
  If \(A\) is a unimodular matrix, then the maximum arithmetic degree of
  \(in_{\vc}(I_A)\) equals the normalized volume of \(\conv(A')\).
  \label{thm:MaxVolume}
\end{thm}

\begin{proof}
  For any \(\vec{c}\), the set of standard pairs of\(in_{\vc}(I_A)\) is
  \(\{(1,\sigma)\ |\ \sigma\ \mbox{is a maximal face of}
  \Delta_{\vc}\}\), and each \((1,\sigma)\) corresponds to the standard
  pair \((1,\sigma\cup\{n+1\})\) of
  \(in_{(\vec{c},0)}(I_{A'})\). Especially, \(\sigma\cup\{n+1\}\) is a
  maximal face of \(\Delta'_{(\vc,0)}\). Therefore,
  \begin{eqnarray*}
    \arithdeg{in_{\vc}(I_A)} &=& \abs{\left\{(1,\sigma)\in
        S\left(in_{\vc}(I_A)\right)\right\}}\\
    &=& \abs{\left\{(1,\sigma\cup\{n+1\})\in
        S\left(in_{(\vc,0)}(I_{A'})\right)\right\}}\\
    &\leq& \abs{\left\{(*,\tau)\in S\left(in_{(\vc,0)}(I_{A'})\right)\
        |\ \tau :\ \mbox{maximal face of}\ \Delta'_{(\vc,0)}\right\}}\\
    &=& \mbox{normalized volume of}\ \conv(A').
  \end{eqnarray*}

  Let \(I_{A}\subset k[\vx]\) and \(I_{A'}\subset
  k[x_1,\ldots,x_n,x_{n+1}]\). Then \(\vx^{\va}-\vx^{\vb}x_{n+1}^k\in
  I_{A'}\) (\(\vx^{\va},\ \vx^{\vb}\in k[\vx]\)) if and only if
  \(\sum_{i=1}^n (a_i-b_i)=k\) and \(\vx^{\va}-\vx^{\vb}\in I_A\). We
  consider that \(\vc=(1,1,\ldots,1)\) and \(\succ\) is any reverse
  lexicographic term order such that \(x_{n+1}\) is the smallest
  variable. Then for any \(g\) in the reduced Gr{\"o}bner basis
  \({\mathcal G}\) for \(I_{A'}\) with respect to \(\succ_{(\vc,0)}\),
  \(in_{\succ_{(\vc,0)}}(g)\) does not contain \(x_{n+1}\) by the
  definition of the term order, and \(in_{\succ_{(\vc,0)}}(g)\) is
  square-free since \(\{in_{\succ_{(\vc,0)}}(g)\ |\ g\in{\mathcal G}\}\)
  minimally generates \(in_{\succ'_{\vc}}(I_A)\) for some term order
  \(\succ'\). Thus the corresponding triangulation
  \(\Delta'_{\succ_{(\vc,0)}}\) is unimodular~\cite{Stu95}, and each
  facet of \(\Delta'_{\succ_{(\vc,0)}}\) corresponds to a standard pair
  of \(in_{\vc}(I_A)\) injectively. Then the arithmetic degree of
  \(in_{\vc}(I_A)\) is equal to the number of facets of
  \(\Delta'_{\succ_{(\vc,0)}}\), which is the normalized volume of
  \(\conv(A')\).
\end{proof}

We consider the primal problem which is equivalent with
\(LP_{A,\vc}(\vb)\):
\[
P_{(M\ I),\widetilde{\vc}}(\widetilde{\vb}) := maximize\
\{(-\widetilde{\vc})\T\vx'\ |\ M\vx'+I_d\vx''=\widetilde{\vb}_{B},\
\vx',\vx''\geq\vzero\},\]
which corresponds to some basis \(B\), and its dual problem
\[
D_{(I\ -M\T),\widetilde{\vb}}(\widetilde{\vc} ) := minimize\
\{\widetilde{\vb}_{B}\T\vy''\ |\ I_{n-d}\vy'-M\T\vy''=\widetilde{\vc},\ 
\vy',\vy''\geq \vzero\},
\]
where \(M\in\ZZ^{d\times(n-d)}\),
\(\widetilde{\vb}=(\widetilde{\vb}_{B},\widetilde{\vb}_{N})\)
(\(\widetilde{\vb}_{B}=(\widetilde{b}_i)_{i\in B}\in\ZZ^d,\
\widetilde{\vb}_{N}=(\widetilde{b}_{i})_{i\notin
  B}=\vzero\in\ZZ^{n-d}\)), \(I_d\in\ZZ^{d\times d}\) and
\(I_{n-d}\in\ZZ^{(n-d)\times(n-d)}\) are identity matrices, \(\vx''\)
(resp. \(\vx'\)) is a basic (resp. non-basic) variable for \(P_{(M\
  I),\widetilde{\vc}}(\widetilde{\vb})\), \(\vy'\) (resp. \(\vy''\)) is
a basic (resp. non-basic) variable for \(D_{(I\
  -M\T),\widetilde{\vb}}(\widetilde{\vc})\), and \(\widetilde{\vc}\) is
a reduced cost for \(B\).

For any standard pair \((1,\sigma)\) of
\(in_{\vc}(I_A)=in_{\widetilde{\vc}}(I_{(M\ I)})\),
\(\overline{\sigma}:=\{1,\ldots,n\}\setminus\sigma\) forms a basis of
the dual problem \(D_{(I\ -M\T),\widetilde{\vb}}(\widetilde{\vc})\)
(Lemma~\ref{lem:Tri_and_Dual}). Let
\[
\sigma_1 := (\{1,\ldots,n\}\setminus B)\cap\sigma,\quad \sigma_2 :=
B\cap\sigma,\quad \overline{\sigma_1} := (\{1,\ldots,n\}\setminus
B)\cap\overline{\sigma},\quad \overline{\sigma_2} :=
B\cap\overline{\sigma}.
\]
Then the {\em reduced\/ cost\/} of \(D_{(I\
  -M\T),\widetilde{\vb}}(\widetilde{\vc})\) for the basis
\(\overline{\sigma}\) is
\[
\widetilde{\vb}'_{\sigma}=\widetilde{\vb}_{\sigma}-N_1\T
(B_1^{-1})\T\widetilde{\vb}_{\overline{\sigma}},\ \mbox{where}\
B_1=(I_{\overline{\sigma_1}}\ (-M\T)_{\overline{\sigma_2}}),\
N_1=(I_{\sigma_1}\ (-M\T)_{\sigma_2}).
\]

\begin{thm}
  The solution of the equation in Algorithm~\ref{algo:IP_Standard_pair}
  (i) for a standard pair \((1,\sigma)\) is the reduced cost of
  \(D_{(I\ -M\T),\widetilde{\vb}}(\widetilde{\vc})\) for the basis
  \(\overline{\sigma}\).
  \label{thm:red_cost}
\end{thm}

\begin{proof}
  We show that \(\widetilde{\vb}'_{\sigma}\) is a solution of the linear 
  system in Algorithm~\ref{algo:IP_Standard_pair} (i) for
  \((1,\sigma)\), i.e. \((M_{\sigma_1}\
  I_{\sigma_2})\widetilde{\vb}'_{\sigma}=\widetilde{\vb}\). This is
  because 
  \begin{eqnarray}
    (M_{\sigma_1}\ I_{\sigma_2})\widetilde{\vb}'_{\sigma} &=&
    (M_{\sigma_1}\ I_{\sigma_2})\widetilde{\vb}_{\sigma} -
    (M_{\sigma_1}\
    I_{\sigma_2})N_1\T(B_1^{-1})\T\widetilde{\vb}_{\overline{\sigma}}
    \nonumber\\
    &=& I_{\sigma_2}\widetilde{\vb}_{\sigma_2} -
    (M_{\sigma_1}(I_{\sigma_1})\T +
    I_{\sigma_2}((-M\T)_{\sigma_2})\T)(B_1^{-1})\T\widetilde{\vb}_{\overline{\sigma}} \nonumber\\
    &=&  I_{\sigma_2}\widetilde{\vb}_{\sigma_2} -
    \left\{(M-M_{\overline{\sigma_1}}(I_{\overline{\sigma_1}})\T) +
      (-M-I_{\overline{\sigma_2}}((-M\T)_{\overline{\sigma_2}})\T)\right\}
    (B_1^{-1})\T\widetilde{\vb}_{\overline{\sigma}} \label{eq:subsec}\\
    &=&  I_{\sigma_2}\widetilde{\vb}_{\sigma_2}
    +(M_{\overline{\sigma_1}}\
    I_{\overline{\sigma_2}})(I_{\overline{\sigma_1}}\
    (-M\T)_{\overline{\sigma_2}})\T
    (B_1^{-1})\T\widetilde{\vb}_{\overline{\sigma}} \nonumber\\
    &=&  I_{\sigma_2}\widetilde{\vb}_{\sigma_2}
    +(M_{\overline{\sigma_1}}\
    I_{\overline{\sigma_2}})B_1\T(B_1^{-1})\T\widetilde{\vb}_{\overline{\sigma}} \nonumber\\
    &=& I_{\sigma_2}\widetilde{\vb}_{\sigma_2} +
    (M_{\overline{\sigma_1}}\ I_{\overline{\sigma_2}})\widetilde{\vb}_{\overline{\sigma}} \nonumber\\
    &=& I_{\sigma_2}\widetilde{\vb}_{\sigma_2} +
    I_{\overline{\sigma_2}}\widetilde{\vb}_{\overline{\sigma_2}}
    \nonumber\\
    &=& \widetilde{\vb},
  \end{eqnarray}
  the equation (\ref{eq:subsec}) follows from the fact that
  \(M=MI=M_{\sigma_1}(I_{\sigma_1})\T+M_{\overline{\sigma_1}}(I_{\overline{\sigma_1}})\T\) and \(-M=I(-M)=I_{\sigma_2}((-M\T)_{\sigma_2})\T+I_{\overline{\sigma_2}}((-M\T)_{\overline{\sigma_2}})\T\).

\end{proof}
  
\section{Gr{\"o}bner bases and standard pairs of the primal minimum cost
  flow problems}
Let \(G_d\) be the acyclic tournament graph with vertices
\(1,2,\ldots,d\) and \(n=\tbinom{d}{2}\) arcs, where each arc \((i,j)\
(i<j)\) is directed from \(i\) to \(j\). We consider the following
minimum cost flow problem \(P_{A,\vc}(\vb)\):
\[
P_{A,\vc}(\vb) := minimize\ \{\vc\T\vx\ |\ A\vx=\vb,\ \vx\geq \vzero\},
\]
where \(A\in\ZZ^{d\times n}\) is the vertex-arc incidence matrix of
\(G_d\).

A {\em walk\/} in \(G_d\) is a sequence \((v_1,v_2,\ldots,v_p)\) of
vertices such that \((v_i,v_{i+1})\) or \((v_{i+1},v_i)\) is an arc of
\(G_d\) for each \(1\leq i<p\). A {\em cycle\/} is a walk
\((v_1,v_2,\ldots,v_p,v_1)\). A {\em circuit\/} is a cycle
\((v_1,v_2,\ldots,v_p,v_1)\) such that \(v_i\neq v_j\) for any \(i\neq
j\).

\begin{defn}
  Let \(C\) be a circuit in \(G_d\) and fix a direction of \(C\). If
  \(C\) passes an arc \((i,j)\) \(u_{ij}^{+}\) times forwardly and
  \(u_{ij}^{-}\) times backwardly, then we define
  \(\vu_{C}^{+}=(u_{ij}^{+})_{1\leq i<j\leq d},
  \vu_{C}^{-}=(u_{ij}^{-})_{1\leq i<j\leq d}\in\RR^n\). The vector
  \(\vu_{C}:=\vu_{C}^{+}-\vu_{C}^{-}\) is called the {\em incidence\/
    vector\/} of \(C\). We identify a cycle \(C\) of \(G_d\) with the
  binomial \(f_{C}:=\vx^{\vu_{C}^{+}}-\vx^{\vu_{C}^{-}}\in I_{A}\).
\end{defn}

\begin{defn}
  A non-zero vector \(\vu\) in \(\ker(A)\) is a {\em circuit\/} if its
  support \(\supp(\vu)\) is minimal with respect to inclusion and the
  elements of \(\vu\) are relatively prime. When \(\vu\in\ker(A)\) is a
  circuit, we also call \(\vx^{\vu^{+}}-\vx^{\vu^{-}}\) a {\em
    circuit\/} of \(I_A\). We denote \({\mathcal C}_A\) a set of
  circuits of \(I_A\).
\end{defn}

Then \({\mathcal C}_{A}\) corresponds to the set of all circuits in
\(G_d\). Let \({\mathcal U}_A\) be the union of all reduced Gr{\"o}bner
bases for \(I_A\) with respect to all term orders, which is called the
{\em universal\/ Gr{\"o}bner\/ basis} of \(I_A\).

\begin{prop}[\cite{Stu95}]
  For the vertex-arc incidence matrix \(A\) of \(G_d\), \({\mathcal
    U}_{A}={\mathcal C}_{A}\). Especially, any reduced Gr{\"o}bner basis
  of \(I_A\) is square-free, and the number of elements in \({\mathcal
    U}_{A}\) is of exponential order with respect to \(d\).
  \label{prop:universalGB}
\end{prop}

\begin{prop}
  \(I_A\) is not homogeneous for the grading \(\deg(x_{i,j})=1\), but is
  homogeneous for the grading \(\deg(x_{i,j})=j-i\).
\end{prop}

\begin{proof}
For any \(d\), \(x_{1,2}x_{2,3}-x_{1,3}\in I_A\) and
\(x_{1,2}x_{2,3}\notin I_A\). This implies that \(I_A\) is not
homogeneous for the grading \(\deg(x_{i,j})=1\).

Let \(v_1,v_2,\ldots,v_p,v_1\) be a circuit in \(G_d\), \(C^{+}:=\{k\ |\ 
v_k<v_{k+1}\}\) and \(C^{-}:=\{k\ |\ v_k>v_{k+1}\}\) (we set
\(v_{p+1}:=v_1\)). The binomial \(f_C\) corresponding to \(C\) is
\(f_C=\prod_{k\in C^{+}}x_{v_kv_{k+1}}-\prod_{k\in
  C^{-}}x_{v_{k+1}v_k}\). Then \(f_C\) is homogeneous for the grading
\(\deg(x_{i,j})=j-i\) because 
\begin{eqnarray*}
  \deg\left(\prod_{k\in
      C^{+}}x_{v_kv_{k+1}}\right)-\deg\left(\prod_{k\in
      C^{-}}x_{v_{k+1}v_k}\right) &=& \sum_{k\in
  C^{+}}(v_{k+1}-v_k)-\sum_{k\in C^{-}}(v_k-v_{k+1})\\
  &=& \sum_{k=1}^{p}(v_{k+1}-v_k)\ =\ 0. 
\end{eqnarray*}
\end{proof}

Thus reduced Gr{\"o}bner basis exists for any
\(\vc\in\RR^n\setminus\{\vzero\}\) by Proposition~\ref{prop:grad}.

\subsection{Some Gr{\"o}bner bases for the primal problem}
We first show that the elements in reduced Gr{\"o}bner bases with
respect to some specific term orders can be given in terms of graphs. As
a corollary, we can show that there exist term orders for which reduced
Gr{\"o}bner bases remain in polynomial order. For other applications of
the Gr{\"o}bner bases found in this section, see our
paper~\cite{IshIma000516}. 

\begin{prop}
  Let \(\succ\) be the purely lexicographic order induced by the
  variable ordering such that \(x_{i,j}\succ x_{k,l}\) if and only if
  \(i<k\) or \((i=k\ \mbox{and}\ j<l)\). Then the reduced Gr{\"o}bner
  basis for \(I_A\) with respect to \(\succ\) is
  \(\{g_{ijk}:=x_{i,j}x_{j,k}-x_{i,k}\ |\  i<j<k\} \cup
  \{g_{ijkl}:=x_{i,k}x_{j,l}-x_{i,l}x_{j,k}\ |\ i<j<k<l\}\). In
  particular, the number of elements in this Gr{\"o}bner basis is equal
  to \(\tbinom{d}{3}+\tbinom{d}{4}\). 
  \label{prop:grob_p_type1}
\end{prop}

\begin{figure}[htbp]
  \begin{center}
    \epsfig{file=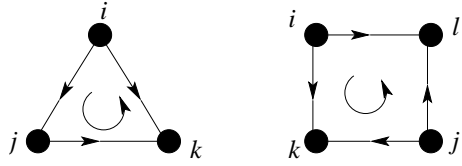,width=6cm}
    \caption{Circuit corresponding to \(g_{ijk}\) (left) and circuit
      corresponding to \(g_{ijkl}\) (right).}
    \label{fig:grob_p_type1}
  \end{center}
\end{figure}

The set \(\{g_{ijk}\ |\  i<j<k\}\) corresponds to all of the circuits of 
length three, and \(\{g_{ijkl}\ |\  i<j<k<l\}\) corresponds to some
circuits of length four uniquely determined for each four vertices
\(i,j,k,l\)~(Figure~\ref{fig:grob_p_type1}).

\begin{proof}
  By Proposition~\ref{prop:universalGB}, it suffices to show that any
  binomial which corresponds to a circuit in \(G_d\) is \(g_{ijk}\),
  \(g_{ijkl}\) or whose initial term is divisible by some
  \(in_{\succ}(g_{ijk})\) or \(in_{\succ}(g_{ijkl})\).

  Any binomial corresponding to a circuit of length \(3\) is contained in
  \(\{g_{ijk}\}\). 
  
  The circuits defined by four vertices \(i<j<k<l\) are
  \(C_1:=(i,j,k,l,i)\), \(C_2:=(i,j,l,k,i)\), \(C_3:=(i,k,j,l,i)\) and
  their opposites. The binomial which corresponds to \(C_1\) or its
  opposite is \(\pm(x_{i,j}x_{j,k}x_{k,l}-x_{i,l})\), whose initial term
  \(x_{i,j}x_{j,k}x_{k,l}\) is divisible by
  \(in_{\succ}(g_{ijk})\). Similarly, the initial term of binomial which
  corresponds to \(C_2\) or its opposite is divisible by
  \(in_{\succ}(g_{ijl})\). The binomial which corresponds to \(C_3\) or
  its opposite is \(g_{ijkl}\).

  Let \(C\) be a circuit of length more than five. Let \(v_1\) be the
  vertex whose label is minimum in \(C\), and
  \(C:=(v_1,v_2,\ldots,v_p,v_1)\). Without loss of generality, we set
  \(v_2<v_p\). Let \(f_C\) be the binomial corresponding to \(C\), then
  \(in_{\succ}(f_C)\) is the product of all variables whose associated
  arcs have the same direction as \((v_1,v_2)\) on \(C\). If
  \(v_2<v_3\), then \((v_1,v_2)\) and \((v_2,v_3)\) have the same
  direction on \(C\). Thus both \(x_{v_1,v_2}\) and \(x_{v_2,v_3}\)
  appear in \(in_{\succ}(f_C)\), and \(in_{\succ}(f_C)\) is
  divisible by \(in_{\succ}(g_{v_1v_2v_3})\). If \(v_2>v_3\), then since
  \(v_3<v_2<v_p\), there exists some \(k\ (3\leq k\leq p-1)\) such that
  \(v_1<v_k<v_2<v_{k+1}\). Then both \(x_{v_1,v_2}\) and
  \(x_{v_k,v_{k+1}}\) appear in \(in_{\succ}(f_C)\), and
  \(in_{\succ}(f_C)\) is divisible by
  \(in_{\succ}(g_{v_1v_kv_2v_{k+1}})\).
\end{proof}

\begin{thm}
  Let \(\succ\) be any term order and
  \(\vc=(c_{1,2},\ldots,c_{1,d},c_{2,3},\ldots,c_{d-1,d})\in\RR^n\)
  satisfy \(c_{i,j}+c_{j,k}>c_{i,k}\) for any \(i<j<k\) and
  \(c_{i,k}+c_{j,l}>c_{i,l}+c_{j,k}\) for any \(i<j<k<l\). Then the
  reduced Gr{\"o}bner basis for \(I_A\) with respect to \(\succ_{\vc}\)
  is the same as the basis in Proposition~\ref{prop:grob_p_type1}.
  \label{thm:grob_p_type1}
\end{thm}

\begin{proof}
  Let \(\succ'\) be the term order defined in
  Proposition~\ref{prop:grob_p_type1}. Then
  \(in_{\succ_{\vc}}(g_{ijk})=x_{i,j}x_{j,k}=in_{\succ'}(g_{ijk})\) since
  \(c_{i,j}+c_{j,k}>c_{i,k}\), and
  \(in_{\succ_{\vc}}(g_{ijkl})=x_{i,k}x_{j,l}=in_{\succ'}(g_{ijkl})\) since
  \(c_{i,k}+c_{j,l}>c_{i,l}+c_{j,k}\). Thus
  \(in_{\succ_{\vc}}(I_A)=in_{\succ'}(I_A)\), which implies that the
  reduced Gr{\"o}bner bases for \(I_A\) with respect to \(\succ_{\vc}\)
  and \(\succ'\) are the same.
\end{proof}

\begin{prop}
  Let \(\succ\) be the purely lexicographic order induced by the
  variable ordering such that \(x_{i,j}\succ x_{k,l}\) if and only if
  \(j-i<l-k\) or \((j-i=l-k\ \mbox{and}\ i<k)\). Then the reduced
  Gr{\"o}bner basis for \(I_A\) with respect to \(\succ\) is
  \(\{g_{ijk}:=x_{i,j}x_{j,k}-x_{i,k}\ |\  i<j<k\} \cup
  \{g_{ijkl}:=x_{i,l}x_{j,k}-x_{i,k}x_{j,l}\ |\ i<j<k<l\}\). In
  particular, the number of elements in this Gr{\"o}bner basis is equal
  to \(\tbinom{d}{3}+\tbinom{d}{4}\).
  \label{prop:grob_p_type2}
\end{prop}

The set \(\{g_{ijk}\ |\  i<j<k\}\) corresponds to all of the circuits of 
length three in \(G_d\), and \(\{g_{ijkl}\ |\ i<j<k<l\}\) corresponds to
the set of circuits of length four in Figure~\ref{fig:grob_p_type1} but
the direction is opposite.

\begin{proof}
  Any binomial corresponds to a circuit of length \(3\) is contained in
  \(\{g_{ijk}\}\). 

  The circuits defined by four vertices \(i<j<k<l\) are
  \(C_1:=(i,j,k,l,i)\), \(C_2:=(i,j,l,k,i)\), \(C_3:=(i,k,j,l,i)\) and
  their opposites. The binomial which corresponds to \(C_1\) or its
  opposite is \(\pm(x_{i,j}x_{j,k}x_{k,l}-x_{i,l})\), whose initial term
  \(x_{i,j}x_{j,k}x_{k,l}\) is divisible by \(in_{\succ}(g_{ijk})\). The
  binomial which corresponds to \(C_2\) or its opposite is
  \(\pm(x_{i,j}x_{j,l}-x_{i,k}x_{k,l})\). If its initial term is
  \(x_{i,j}x_{j,l}\), it is divisible by \(in_{\succ}(g_{ijl})\). On the
  other hand, if initial term is \(x_{i,k}x_{k,l}\), it is divisible by
  \(in_{\succ}(g_{ikl})\). The binomial which corresponds to \(C_3\) or
  its opposite is \(g_{ijkl}\).
  
  Let \(C\) be a circuit of length more than five. Let \((v_1,v_2)\
  (v_1<v_2)\) be the arc in \(C\) which the difference of labels is
  minimum, and \(C:=(v_1,v_2,\ldots,v_p,v_1)\). Let \(f_C\) be the
  binomial corresponding to \(C\), then \(in_{\succ}(f_C)\) is the
  product of all variables whose associated arcs have the same direction
  with \((v_1,v_2)\) on \(C\).
  
  If \(v_2<v_3\), then both \(x_{v_1,v_2}\) and \(x_{v_2,v_3}\) appear
  in \(in_{\succ}(f_C)\), and \(in_{\succ}(f_C)\) is divisible by
  \(in_{\succ}(g_{v_1v_2v_3})\). Similarly, if \(v_p<v_1\), then
  \(in_{\succ}(f_C)\) is divisible by \(in_{\succ}(g_{v_pv_1v_2})\).

  Let \(v_3<v_2\) and \(v_1<v_p\). Then \(v_3<v_1<v_2<v_p\) by the
  definition of \(v_1\) and \(v_2\). If there exists some \(q\) such
  that \(v_q<v_{q+1}<v_{q+2}\), then \(in_{\succ}(f_C)\) is divisible by
  \(in_{\succ}(g_{v_qv_{q+1}v_{q+2}})\). Consider the case that there
  does not exist such \(q\). For any \(s\) such that
  \(v_s<v_1<v_{s+1}<v_2\), \(v_{s+2}<v_{s+1}\) by assumption, and
  \(v_{s+2}<v_1\) by the definition of \(v_1\) and \(v_2\). Thus there
  must be some \(r\ (3\leq r\leq p-1)\) such that
  \(v_r<v_1<v_2<v_{r+1}\) since \(v_3<v_1<v_2<i_p\). Then
  \(in_{\succ}(f_C)\) is divisible by
  \(in_{\succ}(g_{v_rv_1v_2v_{r+1}})\).
\end{proof}

\begin{thm}
  Let \(\succ\) be any term order and
  \(\vc=(c_{1,2},\ldots,c_{1,d},c_{2,3},\ldots,c_{d-1,d})\in\RR^n\) satisfy
  \(c_{i,j}+c_{j,k}>c_{i,k}\) for any \(i<j<k\) and
  \(c_{i,l}+c_{j,k}>c_{i,k}+c_{j,l}\) for any \(i<j<k<l\). Then the
  reduced Gr{\"o}bner basis for \(I_A\) with respect to \(\succ_{\vc}\)
  is the same as the basis in Proposition~\ref{prop:grob_p_type2}.
  \label{thm:grob_p_type2}
\end{thm}

\begin{proof}
  Let \(\succ'\) be the term order defined in
  Proposition~\ref{prop:grob_p_type2}. Then
  \(in_{\succ_{\vc}}(g_{ijk})=x_{i,j}x_{j,k}=in_{\succ'}(g_{ijk})\) since
  \(c_{i,j}+c_{j,k}>c_{i,k}\), and
  \(in_{\succ_{\vc}}(g_{ijkl})=x_{i,l}x_{j,k}=in_{\succ'}(g_{ijkl})\)
  since \(c_{i,l}+c_{j,k}>c_{i,k}+c_{j,l}\). Thus
  \(in_{\succ_{\vc}}(I_A)=in_{\succ'}(I_A)\), which implies that the
  reduced Gr{\"o}bner bases for \(I_A\) with respect to \(\succ_{\vc}\)
  and \(\succ'\) are the same.
\end{proof}

\begin{prop}
  Let \(\succ\) be the purely lexicographic order induced by the
  variable ordering such that \(x_{i,j}\succ x_{k,l}\) if and only if
  \(i<k\) or \((i=k\ \mbox{and}\ j>l)\). Then the reduced Gr{\"o}bner
  basis for \(I_A\) with respect to \(\succ\) is
  \(\{g_{ij}:=x_{i,j}-x_{i,i+1}x_{i+1,i+2}\cdots x_{j-1,j} \ |\
  i<j-1\}\). In particular, the number of elements in this Gr{\"o}bner
  basis is equal to \(\tbinom{d}{2}-(d-1)\).
  \label{prop:grob_p_type3}
\end{prop}

The set \(\{g_{ij}\ |\ i<j-1\}\) corresponds to all of the fundamental
circuits of \(G_d\) for the spanning tree \(T:=\{(i,i+1)\ |\ 1\leq
i<d\}\).

\begin{proof}
  Let \(C\) be a circuit which is not a fundamental circuit of
  \(T\). Let \(v_1\) be the vertex whose label is minimum in \(C\), and
  \(C:=(v_1,v_2,\ldots,v_p,v_1)\). Without loss of generality, we set
  \(v_2<v_p\). Then the variable \(x_{v_1,v_p}\) appears in the initial
  term of the associated binomial \(f_C\). Thus
  \(in_{\succ}(f_C)\) is divisible by \(in_{\succ}(g_{v_1v_p})\).
\end{proof} 

\begin{thm}
  Let \(\succ\) be any term order and
  \(\vc=(c_{1,2},\ldots,c_{1,d},c_{2,3},\ldots,c_{d-1,d})\in\RR^n\)
  satisfy \(c_{i,j}>c_{i,i+1}+c_{i+1,i+2}+\cdots+c_{j-1,j}\) for any
  \(i<j-1\). Then the reduced Gr{\"o}bner basis for \(I_A\) with respect
  to \(\succ_{\vc}\) is the same as the basis in
  Proposition~\ref{prop:grob_p_type3}.
  \label{thm:grob_p_type3}
\end{thm}

\begin{proof}
  Let \(\succ'\) be the term order defined in
  Proposition~\ref{prop:grob_p_type3}. Then
  \(in_{\succ_{\vc}}(g_{ij})=x_{i,j}=in_{\succ'}(g_{ij})\) since
  \(c_{i,j}>c_{i,i+1}+c_{i+1,i+2}+\cdots+c_{j-1,j}\). Thus
  \(in_{\succ_{\vc}}(I_A)=in_{\succ'}(I_A)\), which implies that the
  reduced Gr{\"o}bner bases for \(I_A\) with respect to \(\succ_{\vc}\)
  and \(\succ'\) are the same.
\end{proof}

\subsection{Bounds for the size of Gr{\"o}bner bases}
Generally the degree of any reduced Gr{\"o}bner basis for toric ideal
is of exponential order with respect to the number of rows in the
matrix~\cite{Stu91}, but the cardinality is not well understood. For the
case of the toric ideals of acyclic tournament graphs, since those
vertex-arc incidence matrices are unimodular, the cardinalities of the
reduced Gr{\"o}bner bases may be bounded.

\begin{prop}
  The minimum cardinality of the reduced Gr{\"o}bner bases for \(I_{A}\)
  is \(\tbinom{d}{2}-(d-1)\). The basis we have shown in
  Proposition~\ref{prop:grob_p_type3} is the example achieving this
  cardinality.
\end{prop}

\begin{proof}
  Since the reduced Gr{\"o}bner basis forms a basis for \(I_A\), the
  cardinality of the reduced Gr{\"o}bner basis is more than that of the
  basis for \(I_A\). Since \(I_A\) corresponds to the cycle space of
  \(G_d\), the cardinality of the basis for \(I_A\) is equal to the
  dimension of the cycle space, which is \(\tbinom{d}{2}-(d-1)\).
\end{proof}

To analyze the upper bound for the cardinalities of the reduced
Gr{\"o}bner bases, we calculate all reduced Gr{\"o}bner bases for small
\(d\) using TiGERS~\cite{HubTho99}. Table~\ref{table:result_TiGERS} is
the result for \(d=4,5,6,7\).

\begin{table}[h]
  \begin{center}
    \begin{tabular}{|c||r|r|r|}
      \hline
      {\(d\)} & {\# GB} & {max cardinality} & {min cardinality}\\ \hline
      4 & 10 & 5 & 3 \\ \hline
      5 & 211 & 15 & 6 \\ \hline
      6 & 48312 & 37 & 10 \\ \hline
      7 & \(\geq\) 37665 & \(\geq\) 75 & 15 \\ \hline
    \end{tabular}
  \end{center}
  \caption{Number of reduced Gr{\"o}bner bases, maximum and minimum of
    cardinality.}
  \label{table:result_TiGERS}
\end{table}

For the case of \(d=7\), the number of reduced Gr{\"o}bner bases and the 
maximum of the cardinality are both too large, so we could not know the
correct values.  For \(d\leq 5\), the reduced Gr{\"o}bner basis in
Proposition~\ref{prop:grob_p_type1} is the example achieving maximum
cardinality, but for \(d\geq 6\) the maximum cardinality is a little
larger than the cardinality of Gr{\"o}bner basis in
Proposition~\ref{prop:grob_p_type1}. For \(d=6\), we do not know what
cost vectors produce the Gr{\"o}bner bases of cardinality 37. The
reduced Gr{\"o}bner bases which achieve the maximum cardinality seem to
be complicated and difficult to characterize. 

\begin{prob}
  Are the cardinalities of reduced Gr{\"o}bner bases for \(I_A\)
  of polynomial order with respect to \(d\)?
\end{prob}

\subsection{Standard pairs for primal problem}
\label{subsec:std_pair_primal}
In this section, we assume that \(\vc\) is generic. If \(\vc\) is not
generic, then we may use \(\vc'\) which is obtained by perturbing
\(\vc\) such that \(\vc'\) is generic and \(in_{\vc}(f)\) contains the
term \(in_{\vc'}(f)\) for any \(f\in I_A\). Since one constraint
of \(P_{A,\vc}(\vb)\) is redundant, we can consider the problem
\(P_{\overline{A},\vc}(\overline{\vb})\), which is obtained from
\(P_{A,\vc}(\vb)\) by deleting the last constraint. Then
\(in_{\vc}(I_A)=in_{\vc}(I_{\overline{A}})\), and \(\overline{A}\) is
row-full rank. In addition, the regular triangulation of \(\cone(A)\) and
that of \(\cone(\overline{A})\) by \(\vc\) are the same as a simplicial
complex, thus we denote both triangulations \(\Delta_{\vc}\).

Since any initial ideal \(in_{\vc}(I_A)\) is generated by square-free
monomials (Proposition~\ref{prop:universalGB}), the set of standard
pairs \(S(in_{\vc}(I_A))\) are \((1,\sigma)\) where \(\sigma\) ranges
among all maximal faces of \(\Delta_{\vc}\).

Let \(E\) be a set of arcs in \(G_d\). For \(S\subseteq E\), we
denote \(\vx^{S}:=\prod_{(i,j)\in S}x_{i,j}\). 

The arcs in the optimum flow of uncapacitated minimum cost flow
problem define a forest~\cite{AhuMagOrl93}. Therefore, with the fact the 
dimension of \(\cone(A)\) equals \(d-1\), the next proposition is implied by
Lemma~\ref{lem:Stu_and_STV}, Proposition~\ref{prop:universalGB} and
Lemma~\ref{lem:square_free_std_pairs}.

\begin{prop}
  \((\vx^{\va},\sigma)\) is a standard pair of \(in_{\vc}(I_A)\) if and only 
  if \(\vx^{\va}=1\) and \(\sigma\) is a spanning tree of \(G_d\) such that
  \(\vx^{\sigma}\notin in_{\vc}(I_A)\).
  \label{prop:std_pair_of_digraph}
\end{prop}

Because of the result in Section~\ref{subsec:dual}, there is one-to-one
correspondence between the standard pairs \((1,*)\) of \(in_{\vc}(I_A)\)
and the dual feasible bases of
\(P_{\overline{A},\vc}(\overline{\vb})\). Therefore,
Algorithm~\ref{algo:IP_Standard_pair} for the minimum cost flow problem
\(P_{A,\vc}(\vb)\) is a variant of the enumeration of dual feasible
bases.

Gr{\"o}bner bases which have shown in the previous section give upper
and lower bounds for the arithmetic degree (i.e. bounds for the number
of vertices of the dual polyhedron). The genericity of \(\vc\) implies
that the arithmetic degree of \(in_{\vc}(I_A)\) is equal to or greater
than 1.

\begin{thm}
  The minimum arithmetic degree of \(in_{\vc}(I_A)\) which \(\vc\) varies
  all generic cost vectors equals {\rm 1}.
  \label{thm:min_arith_degree}
\end{thm}

\begin{proof}
  For a cost vector \(\vc\) as in Theorem~\ref{thm:grob_p_type3},
  \(in_{\vc}(I_A)=\langle x_{i,j}\ |\ j-i>1\rangle\). Then
  \(\vx^{\va}\notin in_{\vc}(I_A)\) if and only if \(a_{i,j}=0\) for any
  \((i,j)\) such that \(j-i>1\). The set of all such monomials coincides
  \((1,\{(1,2),(2,3),\ldots,(d-1,d)\})\). Thus only this pair is a
  standard pair of \(in_{\vc}(I_A)\). 
\end{proof}

\setcounter{section}{4}
\setcounter{thm}{14}
\begin{thm}
  The maximum arithmetic degree of \(in_{\vc}(I_A)\) which \(\vc\) varies
  all generic cost vectors equals
  \(C_{d-1}:=\frac{1}{d}\tbinom{2(d-1)}{d-1},\) which is the \((d-1)\)-th
  Catalan number.
  \label{thm:max_arith_degree}
\end{thm}

To show this theorem, we use the next result due to Gelfand et
al.~\cite{GelGraPos96} which studies about some hypergeometric
function.

\begin{lem}[\cite{GelGraPos96}]
  Let \(A'\) be the enlarged matrix {\rm (\ref{eq:homog_matrix})} for the
  incidence matrix \(A\) of the acyclic tournament graph with \(d\)
  vertices, and \(\conv(A')\) be the convex hull of
  \(\va'_1,\ldots,\va'_{n+1}\). Then the normalized volume of \(\conv(A')\)
  equals the \((d-1)\)-th Catalan number \(C_{d-1}\).
  \label{lem:normal_vol_of_homog}
\end{lem}

\begin{proof2}{Theorem~\ref{thm:max_arith_degree}}
  Since \(A\) is unimodular, \(\arithdeg{in_{\vc}(I_A)}\leq
  (\mbox{normalized volume of}\ \conv(A'))\) \(= C_{d-1}\) by
  Theorem~\ref{thm:MaxVolume}. 

  Because of Proposition~\ref{prop:std_pair_of_digraph} and
  Theorem~\ref{thm:grob_p_type1}, for \(\vc\) as in
  Theorem~\ref{thm:grob_p_type1}, \((1,\sigma)\) is a standard pair of
  \(in_{\vc}(I_A)\) if and only if \(\sigma\) is a spanning tree of the
  acyclic tournament graph which satisfies the following two conditions: 
  \begin{description}
   \setlength{\itemsep}{0pt}
  \item[(a)] there are no \(1\leq i<j<k\leq d\) such that both \((i,j)\)
    and \((j,k)\) are arcs in \(\sigma\), and
  \item[(b)] there are no \(1\leq i<j<k<l\leq d\) such that both \((i,k)\)
    and \((j,l)\) are arcs in \(\sigma\).
  \end{description}
  The number of such spanning trees are known to be the \((d-1)\)-th
  Catalan number (e.g. see \cite{Sta99Vol2}).
\end{proof2}
    
We remark that the Catalan number equals
\(\frac{4^n}{\sqrt{\pi}n^{3/2}}\left(1+O\left(\frac{1}{n}\right)\right)\)
(e.g. see \cite{CorLeiRiv90}).


\section{Gr{\"o}bner bases and standard pairs of dual minimum cost
  flow problems}
In this section, we analyze Gr{\"o}bner bases and standard pairs for
the dual minimum cost flow problems. As in Section~\ref{subsec:dual}, we
study the problem which is equivalent with \(P_{A,\vc}(\vb)\):
\[
P_{(M\ I),\widetilde{\vc}}(\widetilde{\vb}) := maximize\
\{(-\widetilde{\vc})\T\vx'\ |\ M\vx'+I\vx''=\widetilde{\vb}_B,\
\vx',\vx''\geq\vzero\},
\]
which corresponds to the basis \(\{(1,2),(2,3),\ldots,(d-1,d)\}\), and its 
dual problem
\[
D_{(I\ -M\T),\widetilde{\vb}}(\widetilde{\vc} ) := minimize\
\{\widetilde{\vb}_B\T\vy''\ |\ I\vy'-M\T\vy''=\widetilde{\vc},\
\vy',\vy''\geq\vzero\},
\]
where \((M\ I)\) (resp. \((I\ -M\T)\)) is the fundamental cutset
(resp. fundamental circuit) matrix which corresponds to the spanning
tree \(\{(1,2),(2,3),\ldots,(d-1,d)\}\), \(\widetilde{\vc}\) is the
reduced cost corresponding to the basis
\(\{(1,2),(2,3),\ldots,(d-1,d)\}\),
\(\widetilde{\vb}=(\widetilde{\vb}_B,\widetilde{\vb}_N)=(\widetilde{b}_{ij})_{1\leq
  i<j\leq d}\),\ \(\widetilde{\vb}_B=(\widetilde{b}_{i,i+1})_{1\leq
  i<d}\),\ \(\widetilde{\vb}_{N}=(\widetilde{b}_{i,j})_{i<j-1}=\vzero\), 
and
\[
 \begin{array}{l}
  \vx=(\vx',\vx''),\ \vx'=(x_{1,3}\ldots,x_{1,d},x_{2,4},\ldots,x_{d-2,d}),\
  \vx''=(x_{1,2},x_{2,3},\ldots,x_{d-1,d})\\
  \vy=(\vy'',\vy'),\ \vy'=(y_{1,3}\ldots,y_{1,d},y_{2,4},\ldots,y_{d-2,d}),\ 
  \vy''=(y_{1,2},y_{2,3},\ldots,y_{d-1,d}).
\end{array}
\]
Then \(P_{(M\ I),\widetilde{\vc}}(\widetilde{\vb})\) has \(d-1\)
constraints~(i.e. \((M\ I)\in\ZZ^{(d-1)\times n}\)), \(D_{(I\
  -M\T),\widetilde{\vb}}(\widetilde{\vc})\) has \(n-d+1\)
constraints~(i.e. \((I\ -M\T)\in\ZZ^{(n-d+1)\times n}\)).

Let \(G_d=(V,E)\). \(D\subseteq E\) is a {\em cutset\/} in \(G_d\) if there
exists a partition \((V_1,V_2)\) of \(V\)\ (i.e. \(V_1\cap
V_2=\emptyset,\ V_1\cup V_2=V\)) such that \(D=\{(i,j)\in E\ |\ i\in
V_1\ \mbox{and}\ j\in V_2,\ \mbox{or}\ i\in V_2\ \mbox{and}\ j\in V_1\}\).

\begin{defn}
  Let \(D\) be a cutset in \(G_d\) which corresponds to \(V=(V^{+},V^{-})\). We
 define the vector
  \(\vu_{D}\in\RR^n\) as
  \[
   (\vu_{D})_{ij}:=\left\{
    \begin{array}{ll}
      1 & (i\in V^{+}\ \mbox{and}\ j\in V^{-})\\
      -1 & (i\in V^{-}\ \mbox{and}\ j\in V^{+})\\
      0 & (\mbox{otherwise})
    \end{array}\right. .\]
  The vector \(\vu_{D}\) is called the {\em incidence\/ vector\/} of \(D\).
\end{defn}

We identify a cutset \(D\) which corresponds to \((V^{+},V^{-})\) with
the binomial \(f_{D}:=\vx^{\vu_{D}^{+}}-\vx^{\vu_{D}^{-}}\). Since the
rank of the fundamental circuit matrix \((I\ -M\T)\) is \(n-d+1\) and
each row vector of the fundamental cutset matrix \((M\ I)\) is in
\(\ker((I\ -M\T))\), the set of row vectors of the fundamental cutset
matrix \((M\ I)\) forms a basis of \(\ker(I\ -M\T)\). 

For the fundamental circuit matrix \((I\ -M\T)\), the set of circuits
\({\mathcal C}_{(I\ -M\T)}\) corresponds to the set of all cutsets of
\(G_d\). Since the fundamental circuit matrix \((I\ -M\T)\) is totally
unimodular (e.g. see~\cite{TakFuj81}),
Proposition~\ref{prop:universalGB} implies \({\mathcal C}_{(I\
  -M\T)}={\mathcal U}_{(I\ -M\T)}\).

\begin{prop}
  For a cost vector \(\widetilde{\vb}\) such that the linear system 
  \((M\ I)\vx=\widetilde{\vb}_B\) has a non-negative solution, \(I_{(I\
    -M\T)}\) has a reduced Gr{\"o}bner basis with respect to
  \(\widetilde{\vb}\).
  \label{prop:dual_cost}
\end{prop}

\begin{proof}
  Let \(\va\geq\vzero\) be a solution of \((M\
  I)\vx=\widetilde{\vb}_B\). We denote \(\vr_i\) the \(i\)-th row of
  \((M\ I)\), i.e. the row which corresponds to the fundamental cutset
  for the arc \((i,i+1)\). For each cutset \(D\) corresponds to
  \((V^{+},V\setminus V^{+})\)\ (\(V^{+}\subseteq \{1,\ldots,d-1\}\)),
  since \(\vu_{D}=\sum_{i\in V^{+},\ i+1\notin V^{+}}\vr_i-\sum_{i\notin
    V^{+},\ i+1\in V^{+}}\vr_i\), 
  \begin{eqnarray*}
    \va\cdot\vu_{D} &=& \sum_{i\in V^{+},\ i+1\notin
      V^{+}}\va\cdot\vr_i-\sum_{i\notin V^{+},\ i+1\in
      V^{+}}\va\cdot\vr_i\\
    &=& \sum_{i\in V^{+},\ i+1\notin
      V^{+}}\widetilde{b}_{i,i+1}-\sum_{i\notin V^{+},\ i+1\in
      V^{+}}\widetilde{b}_{i,i+1}\\
    &=& \widetilde{\vb}\cdot\vu_{D}.
  \end{eqnarray*}
  Thus \(in_{\va}(f_{D})=in_{\widetilde{\vb}}(f_{D})\) for any cutset
  \(D\), and
  \(in_{\va}(I_{(I\ -M\T)})=in_{\widetilde{\vb}}(I_{(I\ -M\T)})\). Since
  \(\va\geq\vzero\), \(I_{(I\ -M\T)}\) has a reduced Gr{\"o}bner basis with
  respect to \(\widetilde{\vb}\).
\end{proof}

\setcounter{section}{2}
\setcounter{thm}{1}
\begin{exmp}[continued]
  Let \(\vc=(3,1,2)\) and \(\vb=(4,5,-9)\). Then the primal and dual
  problem which corresponds to the spanning tree \(\{(1,2),(2,3)\}\) are
  the following.
  \[
  \begin{array}{lllll}
    \max & 4x_{1,3} &\qquad& \min & 4y_{1,2}+9y_{2,3}\\
    \hbox{s.t.} & \left(
      \begin{array}{c|cc} 1 & 1 & 0\\ 1 & 0 & 1\end{array}
    \right)\left(
      \begin{array}{c} x_{1,3} \\ x_{1,2} \\ x_{2,3} \end{array}
    \right)=\left(\begin{array}{c} 4\\ 9\end{array}\right) & &
    \hbox{s.t.} & \left(
      \begin{array}{c|cc} 1 & -1 & -1\end{array}\right)\left(
      \begin{array}{c} y_{1,3} \\ y_{1,2} \\ y_{2,3}
      \end{array}\right)=-4\\
    & x_{1,2},x_{1,3},x_{2,3}\ge 0 & & & y_{1,2},y_{1,3},y_{2,3}\ge 0
  \end{array}
  \]
  Then \(I_{(1,-1,-1)}=\langle x_{1,2}-x_{2,3},\ x_{1,2}x_{1,3}-1,\
  x_{1,3}x_{2,3}-1\rangle\) and reduced Gr{\"o}bner basis for
  \(\widetilde{\vb}=(4,0,9)\) is \(\{x_{2,3}-x_{1,2},\
  x_{1,2}x_{1,3}-1\}\).
\end{exmp}

\setcounter{section}{5}
\subsection{Gr{\"o}bner basis for dual problems}
As for primal problems, we show that the elements in reduced Gr{\"o}bner 
basis to some specific term order can be given in terms of graphs.

\setcounter{thm}{2}
\begin{thm}
  Let \(\widetilde{\vb}\) be the cost vector which satisfies the
  condition in Proposition~\ref{prop:dual_cost},
  \(\widetilde{b}_{i,i+1}>\widetilde{b}_{j,j+1}\) \((1\leq\any i<\any
  j\leq d)\) and \(\widetilde{b}_{i,j}=0\) {\rm (\(\any i, j\) such that
    \(j>i+1\))}. Then the reduced Gr{\"o}bner basis for \(I_{(I\
    -M\T)}\) with respect to \(\widetilde{\vb}\) is
  \(\left\{g_{i}:=\prod_{j<i}x_{j,i}-\prod_{k>i} x_{i,k}\ |\
    i=2,3,\ldots,d\right\}\). In particular, the number of elements in
  this Gr{\"o}bner basis is equal to \(d-1\).
  \label{thm:grob_d}
\end{thm}

Each \(g_{i}\) is an incidence vector of the cutset which corresponds to 
\((V\setminus\{i\},\{i\})\). 

\begin{proof}
  For a cutset \(D\) which corresponds to \((V^{+},V^{-})\) such that
  \(1\in V^{+}\), we define \(P^{+}:=\{i\in V^{+}\ |\ i\neq d,\ i+1\in
  V^{-}\}\) and \(P^{-}:=\{i\in V^{-}\ |\ i\neq d,\ i+1\in
  V^{+}\}\). Let \(P^{+}=\{i_1,\ldots,i_p\}\)\ (\(i_1<i_2<\cdots<i_p\))
  and \(P^{-}=\{j_1,\ldots,j_q\}\)\ (\(j_1<j_2<\cdots<j_q\)). Then
  \(p=q\) or \(p=q+1\), and
  \(i_1<j_1<i_2<j_2<\cdots<i_k<j_k<i_{k+1}<j_{k+1}<\cdots\). Since
  \(\widetilde{\vb}\cdot\vu_{D}^{+} = \sum_{r=1}^{p}
  \widetilde{b}_{i_r,i_r+1} > \sum_{r=1}^{q} \widetilde{b}_{j_r,j_r+1} =
  \widetilde{\vb}\cdot\vu_{D}^{-}\),
  \(in_{\widetilde{\vb}}(f_{D})=\vx^{\vu_{D}^{+}}\). Since
  \(in_{\widetilde{\vb}}(g_{i_1+1})=\prod_{j\leq i_1} x_{j,i_1+1}\),
  \(in_{\widetilde{\vb}}(f_{D})\) can be reduced by
  \(in_{\widetilde{\vb}}(g_{i_1+1})\).
\end{proof}

\subsection{Bounds for the size of Gr{\"o}bner bases}
\begin{prop}
  The minimum cardinality of the reduced Gr{\"o}bner bases for
  \(I_{(I\ -M\T)}\) is \(d-1\). The basis we have shown in
  Theorem~\ref{thm:grob_d} is the example achieving this cardinality.
\end{prop}

\begin{proof}
  Since the reduced Gr{\"o}bner basis forms a basis for \(I_{(I\ -M\T)}\),
  the cardinality of the reduced Gr{\"o}bner basis is more than that of
  the basis for \(I_{(I\ -M\T)}\), which is \(d-1\). 
\end{proof}

To analyze the upper bound for the cardinalities of the reduced
Gr{\"o}bner bases, we calculate all reduced Gr{\"o}bner bases for small
\(d\) using TiGERS~\cite{HubTho99}. Table~\ref{table:result_TiGERS_dual} 
is the result for \(d=4,5,6,7\).

\begin{table}[h]
  \begin{center}
    \begin{tabular}{|c||r|r|r|}
      \hline
      {\(d\)} & {\# GB} & {max cardinality} & {min cardinality}\\ \hline
      4 & 7 & 5 & 3 \\ \hline
      5 & 48 & 10 & 4 \\ \hline
      6 & 820 & 20 & 5 \\ \hline
      7 & 44288 & 39 & 6\\ \hline
    \end{tabular}
  \end{center}
  \caption{Number of reduced Gr{\"o}bner bases of dual problems, maximum
    and minimum of cardinality.}
  \label{table:result_TiGERS_dual}
\end{table}

We do not know what cost vectors produce the Gr{\"o}bner bases of
maximum cardinality. The reduced Gr{\"o}bner bases which achieve the maximum 
cardinality seem to be complicated and difficult to characterize. 

\begin{prob}
  Are the cardinalities of reduced Gr{\"o}bner bases for \(I_{(I\ -M\T)}\)
  of polynomial order with respect to \(d\)?
\end{prob}

\subsection{Standard pairs for dual problem}
In this section, we assume that \(\widetilde{\vb}\) is generic same
as Section~\ref{subsec:std_pair_primal}. Since any initial ideal
\(in_{\widetilde{\vb}}(I_{(I\ -M\T)})\) is generated by square-free
monomials, any standard pair in \(S(in_{\widetilde{\vb}}(I_{(I\
  -M\T)}))\) is of the form \((1,*)\). Moreover, since the support of
each optimal solution of \(D_{(I\
  -M\T),\widetilde{\vc}}(\widetilde{\vb})\) does not include a cutset,
with the fact that \(\dim\cone((I\ -M\T))=n-d+1\), the next proposition
is implied by Lemma~\ref{lem:Stu_and_STV} and
Lemma~\ref{lem:square_free_std_pairs}.

\begin{prop}
  \((\vx^{\va},\sigma)\) is a standard pair of
  \(in_{\widetilde{\vb}}(I_{(I\ -M\T)})\) if and only if \(\vx^{\va}=1\)
  and \(\sigma\) is a co-tree of \(G_d\) such that \(\vx^{\sigma}\notin
  in_{\widetilde{\vb}}(I_{(I\ -M\T)})\).
  \label{prop:std_pair_of_dual_digraph}
\end{prop}

\setcounter{section}{2}
\setcounter{thm}{1}
\begin{exmp}[continued]
  For \(\vc=(3,1,2)\) and \(\vb=(4,5,-9)\), the initial ideal
  \(in_{(4,0,9)}(I_{(1|-1,-1)}) = \langle x_{2,3},\
  x_{1,2}x_{1,3}\rangle\) has two standard pairs \((1,\{(1,2)\})\) and
  \((1,\{(1,3)\})\).
\end{exmp}

\setcounter{section}{5}
\setcounter{thm}{6}
\begin{thm}
  For any \(\widetilde{\vb}\) which satisfies the
  condition in Proposition~\ref{prop:dual_cost}, there exists \(S\subset
  \{1,\ldots,d-1\}\) with \(\abs{S}\geq\lfloor(d-1)/6\rfloor\) such
  that, for any \(\sigma\subseteq S\), there exists a spanning tree
  \(T_{\sigma}\) of \(G_d\) which satisfies the following:
  \begin{description}
    \setlength{\itemsep}{0pt}
  \item[(A)] \(T_{\sigma}\) contains the arc set \(\{(i,i+1)\ |\ i\in
    S\setminus\sigma\}\) and does not contain any arc in \(\{(j,j+1)\ |\
    j\in\sigma\}\), and
  \item[(B)] \((1,\overline{T_{\sigma}})\) is a standard pair of
    \(in_{\widetilde{\vb}}(I_{(I\ -M\T)})\), where
    \(\overline{T_{\sigma}}:=E\setminus T_{\sigma}\) is a co-tree of
    \(T_{\sigma}\).
  \end{description}
  Especially, since \(T_{\sigma}\neq T_{\tau}\) for any
  \(\sigma,\tau\subseteq S\)\ \((\sigma\neq\tau)\),
  \(in_{\widetilde{\vb}}(I_{(I\ -M\T)})\) has at least
  \(\Omega(2^{\lfloor d/6\rfloor})\) standard pairs for any generic
  \(\vb\) which satisfies the condition in
  Proposition~\ref{prop:dual_cost}. 
  \label{thm:std_pair_dual}
\end{thm}

\begin{proof}
We divide \(\{1,\ldots,d-1\}\) into the following subsets.
\[
\begin{array}{l}
M_0:=\{i\in\{1,\ldots,d-1\}\ |\ x_{i,i+1}\in in_{\widetilde{\vb}}(I_{(I\
  -M\T)})\}\\
M_1:=\{i\in\{1,\ldots,d-1\}\ |\ i\notin M_0,\ i\equiv 0\ ({\rm mod}\
3)\}\\
M_2:=\{i\in\{1,\ldots,d-1\}\ |\ i\notin M_0,\ i\equiv 1\ ({\rm mod}\
3)\}\\
M_3:=\{i\in\{1,\ldots,d-1\}\ |\ i\notin M_0,\ i\equiv 2\ ({\rm mod}\
3)\}\\
\end{array}
\]

\begin{lem}
  \(\abs{M_0}\leq \lceil (d-1)/2\rceil.\)
  \label{lem:abs_M0}
\end{lem}

\begin{proof2}{Lemma~\ref{lem:abs_M0}}
  We consider a cutset \(D\) which corresponds to \((V^{+},V^{-})\) such 
  that \(f_D\) contains \(x_{i,j}\) as a term of degree \(1\). Without
  loss of generality, we set \(i\in V^{+}\). We assume that
  \(j-i>1\). Then for any \(k\ (i<k<j)\), if \(k\in V^{+}\) then \(f_D\)
  contains \(x_{k,j}\) and \(x_{i,j}\) in the same term, otherwise
  \(f_D\) contains \(x_{i,k}\) and \(x_{i,j}\) in the same term, which
  contradicts that \(x_{i,j}\) is a term of \(f_D\) of degree 1. Thus
  \(j=i+1\). In addition, \(k\in V^{-}\) for any \(k<i\) and \(k\in
  V^{+}\) for any \(k>i+1\). Therefore, \(V^{+}=\{i,i+2,i+3,\ldots,d\}\)
  and \(V^{-}=\{1,\ldots,i-1,i+1\}\).

  We consider that \(in_{\widetilde{\vb}}(f)=x_{i,i+1}\) for some
  \(f\in I_{(I\ -M\T)}\). If \(x_{i-1,i}\in in_{\widetilde{\vb}}(I_{(I\
    -M\T)})\), then \(f\) can be reduced by the binomial corresponding
  to the cutset between \(\{i-1,i+1,\ldots,d\}\) and
  \(\{1,\ldots,i-2,i\}\) to 
  \[
  f':=x_{i,i+1}-\left(\prod_{k\leq
      i-2}x_{k,i}\right)\left(\prod_{k\geq
      i+2}x_{i+1,k}\right)\left(\prod_{\nidan{k\leq i-1,}{l\geq
        i+2}}x_{k,l}\right)\left(\prod_{k\leq
      i-2}x_{k,i-1}\right)\left(\prod_{k\geq
      i+1}x_{i,k}\right)\left(\prod_{\nidan{k\leq i-2,}{l\geq
        i+1}}x_{k,l}\right),
  \]
  and its initial term is \(x_{i,i+1}\). Since both terms of this
  binomial contain \(x_{i,i+1}\), this implies that
  \(in_{\widetilde{\vb}}(f'/x_{i,i+1})=1\). Since \(\widetilde{\vb}\)
  defines a term order by Proposition~\ref{prop:dual_cost}, this is a 
  contradiction.

  Similarly, \(x_{i+1,i+2}\notin
  in_{\widetilde{\vb}}(I_{(I\ -M\T)})\). Thus \(\abs{M_0}\leq \lceil
  (d-1)/2\rceil\).
\end{proof2}

Thus at least one of \(M_1,\ M_2,\ M_3\) has at least \(\lfloor
(d-1)/6\rfloor\) elements. Let \(S\) be one of such \(M_i\
(i=1,2,3)\). For any \(\sigma:=\{i_1>i_2>\cdots>i_r\}\subseteq S\), we
construct desired spanning trees \(T_{\emptyset},\ T_{\{i_1\}},
T_{\{i_1,i_2\}},\ \ldots,\ T_{\sigma}\) inductively. 

\medskip
\noindent
{\bf \(\cdot\) Initial step:}

Let \(T_{\emptyset}:=\{(1,2),(2,3),\ldots,(d-1,d)\}\). Clearly
\(T_{\emptyset}\) is a spanning tree. Since the reduced Gr{\"o}bner
basis corresponds to a subset of cutsets, the initial term of any
elements of reduced Gr{\"o}bner basis contains a variable \(x_{i,i+1}\)
for some \(i\). Thus \(\vx^{\overline{T_{\emptyset}}}\notin
in_{\widetilde{\vb}}(I_{(I\ -M\T)})\).

\medskip
\noindent
{\bf \(\cdot\) Induction step:}

Let \(T_{\sigma\setminus\{i_r\}}\) be the desired spanning tree for
\(\sigma\setminus\{i_r\}\). We define two edge set
\[
T^{1} :=
\left\{T_{\sigma\setminus\{i_r\}}\setminus\left\{(i_r,i_r+1)\right\}\right\}
\cup\left\{(i_r,i_r+2)\right\},\quad
T^{2}:=\left\{T^{1}\setminus\left\{(i_r+1,i_r+2)\right\}\right\}\cup\left\{(i_r-1,i_r+1)\right\}.
\]
Then both \(T^{1}\) and \(T^{2}\) are spanning trees and satisfy the
condition (A). We show that either \(T^{1}\) or \(T^{2}\) satisfies the
condition (B).

\noindent
\underline{(a)\ The case that \(T^{1}\) satisfies the condition (B).}\
Then \(T^{1}\) is a desired spanning tree \(T_{\sigma}\).

\noindent
\underline{(b)\ The case that \(T^{1}\) does not satisfy the condition
  (B).}\\
In this case \(\vx^{\overline{T^{1}}}\in
in_{\widetilde{\vb}}(I_{(I\ -M\T)})\). Let \({\mathcal G}\) be the
reduced Gr{\"o}bner basis for \(I_{(I\ -M\T)}\) with respect to
\(\widetilde{\vb}\). Then \(\vx^{\overline{T^{1}}}\) can be reduced some 
binomial \(g\in{\mathcal G}\), and such \(g\) is one of the following
form (See Figure~\ref{Fig:dual_std_pair_1}).
\begin{description}
  \setlength{\itemsep}{0pt}
\item[(i)] \(g^{(1)}_{(p)}\) which corresponds to the cutset for
  \((V^{+},V^{-})\),
  \(V^{+}=\{p,p+1,\ldots,i_r,i_r+2,i_r+3,\ldots,d\}\) and
  \(V^{-}=\{1,2,\ldots,p-1,i_r+1\}\) for some \(p\leq i_r\), and its
  initial term is a product of variables corresponds to arcs from
  \(V^{+}\) to \(V^{-}\), or
\item[(ii)\ (The case of \(r>1\))] \(g^{(2)}_{(p,t)}\) which corresponds
  to the cutset for \((V^{+},V^{-})\),
  \(V^{-}=\{1,2,\ldots,p-1,i_r+1,i_{q(1)}+1,\ldots,i_{q(t)}+1\}\) and
  \(V^{+}=V\setminus V^{-}\) for \(1\leq\fitter q(t)<\cdots<\fitter
  q(1)<r\) such that \((i_{q(k)}+1,i_{q(k)}+2)\in
  T_{\sigma\setminus\{i_r\}}\) for \(k=1,\ldots,t\) and \(1\leq\fitter
  p\leq i_r\), and its initial term is a product of variables
  corresponds to arcs from \(V^{+}\) to \(V^{-}\).
\end{description}
\begin{figure}[h]
  \begin{center}
    \epsfig{file=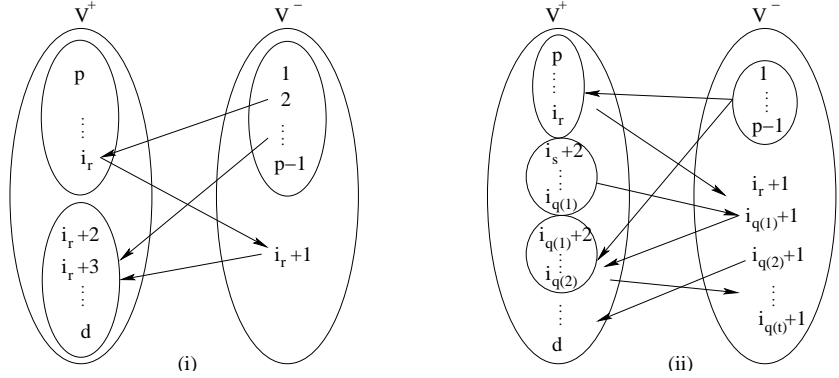,width=11cm}
    \caption{Cutsets corresponding to binomials (i) and (ii).}
    \label{Fig:dual_std_pair_1}
  \end{center}
\end{figure}

\vspace*{-0.5cm}
\begin{lem}
  \(g^{(1)}_{(p)}\in{\mathcal G}\) for some \(p\) and
  \(\vx^{\overline{T^{1}}}\) can be reduced by \(g^{(1)}_{(1)}\),
  i.e. the initial term of \(g^{(1)}_{(1)}\) corresponds to the set of
  arcs \(\{(k,i_r+1) : k\leq i_r\}\).
  \label{lem:dual_std_pair_1}
\end{lem}

\begin{proof2}{Lemma~\ref{lem:dual_std_pair_1}}
  The case of \(r=1\) is trivial.
  
  We suppose \(r>1\) and \(\vx^{\overline{T^{1}}}\) cannot be reduced by
  any \(g^{(1)}_{(p)}\). Then \(\vx^{\overline{T^{1}}}\) can be reduced
  by some \(g^{(2)}_{(p,t)}\) which is an element of \({\mathcal G}\),
  and \(\vx^{\overline{T^{1}}}\) can be also reduced by
  \(g^{(2)}_{(1,t)}\) (otherwise, \(g^{(2)}_{(p,t)}\) is reduced by
  \(g^{(2)}_{(1,t)}\) and \(g^{(2)}_{(p,t)}\) cannot be an element in
  \({\mathcal G}\)).
  
  Suppose that \(\vx^{\overline{T^{1}}}\) can be reduced by
  \(g^{(2)}_{(1,t)}\) with \(t=1\). Let \(m_1\) be the monomial obtained
  by reducing \(\vx^{\overline{T^{1}}}\) by \(g^{(2)}_{(1,t)}\), then
  \(m_1\) can be reduced to the monomial \(m_2\) by \(g^{(1)}_{(1)}\)
  (the initial term of \(g^{(1)}_{(1)}\) is a product of variables
  corresponds to arcs from \(V^{-}\) to \(V^{+}\) by assumption).

  \begin{table}[h]
    \begin{center}
      \begin{tabular}{|l|l|l|l|}
        \hline
        \multicolumn{2}{|c|}{\em reduce\/ by\/ \(g^{(2)}_{(1,1)}\)\/} &
        \multicolumn{2}{c|}{\em reduce\/ by\/ \(g^{(1)}_{(1)}\)\/}\\
        \hline
        divided variables & multiplied variables & divided variables &
        multiplied variables\\
        \hline
        \(\{x_{k,i_r+1}: k\leq i_r\},\) & \(\{x_{i_r+1,l} :\) &
        \(\{x_{i_r+1,l} :\) & \(\{x_{k,i_r+1}: k\leq i_r\}\)\\
        \(\{x_{k,i_{q(1)}+1} :\) & \(\quad l\geq i_r+2,\ l\neq
        i_{q(1)}+1\},\) & \(\quad l\geq i_r+2\}\) & \\
        \(\quad k\leq i_{q(1)},\ k\neq i_r+1\}\) & \(\{x_{i_{q(1)}+1,l}: 
         l\geq i_{q(1)}+2\}\) & & \\
        \hline
      \end{tabular}
    \end{center}
    \label{table:del_in_vars}
    \caption{Divided and multiplied variables while reducing by
      \(g^{(2)}_{(1,1)}\) and \(g^{(1)}_{(1)}\).}
  \end{table}

  For a binomial \(f_{D}\in I_{(I\ -M\T)}\) which corresponds to the
  cutset \(D\) for \((V_{D}^{+},V_{D}^{-})\) such that
  \(V_{D}^{-}=\{i_{q(1)}+1\}\) and \(V_{D}^{+}=V\setminus V_{D}^{-}\),
  \(in_{\widetilde{\vb}}(f_{D})\) corresponds to arcs from \(V_{D}^{-}\)
  to \(V_{D}^{+}\) (otherwise,
  \(\vx^{\overline{T_{\sigma\setminus\{i_r\}}}}\) can be reduced by
  \(f_D\), which contradicts the assumption of the induction).  Then
  \(m_2\) can be reduced by \(f_D\), and the resulting monomial is
  \(\vx^{\overline{T^{1}}}\)~(see Table 3), which contradicts to the
  definition of a term order by \(\widetilde{\vb}\).

  Similarly, in the case that \(\vx^{\overline{T^{1}}}\) can be reduced
  by \(g^{(2)}_{(1,t)}\) for some \(t>1\), using \(f_D\in I_{(I\ -M\T)}\)
  which corresponds to the cutset \(D\) for \((V_{D}^{+},V_{D}^{-})\)
  such that \(V_{D}^{-}=\{i_{q(1)}+1,i_{q(2)}+1,\ldots,i_{q(t)}+1\}\),
  and \(V_{D}^{+}=V\setminus V_{D}^{-}\), we can show a
  contradiction. Thus there exists some \(p\) such that
  \(g^{(1)}_{(p)}\in{\mathcal G}\).

  If \(\vx^{\overline{T^{1}}}\) cannot be reduced by \(g^{(1)}_{(1)}\),
  i.e. the initial term of \(g^{(1)}_{(1)}\) corresponds to the set of
  arcs \(\{(i_r+1, l):\ l\geq i_r+2\}\), then \(g^{(1)}_{(p)}\) can be
  reduced by \(g^{(1)}_{(1)}\), which contradicts that \(g^{(1)}_{(p)}\)
  is an element of reduced Gr{\"o}bner basis \({\mathcal G}\). Thus the
  second statement follows.
\end{proof2}

\medskip
We show that if \(\vx^{\overline{T^{1}}}\in
in_{\widetilde{\vb}}(I_{(I\ -M\T)})\), then \(\vx^{\overline{T^{2}}}\)
cannot be reduced by any binomial in \({\mathcal G}\). If
\(\vx^{\overline{T^{2}}}\) can be reduced by some \(g\in {\mathcal G}\),
then such \(g\) is one of the following form.
\begin{description}
  \setlength{\itemsep}{0pt}
\item[(i)] the binomial \(g^{(1)}_{(i_r)}\), and its initial term is
  \(x_{i_r,i_r+1}\),
\item[(ii)] any binomial which corresponds to the cutset for 
  \((V^{+},V^{-})\) such that \(i_r+1\in V^{+}\) and
  \(1,2,\ldots,i_r,\) \(i_r+2\in V^{-}\), and its initial term is a product
  of variables correspond to arcs from \(V^{+}\) to \(V^{-}\), or
\item[(iii)\ (The case of \(r>1\))] \(g^{(2)}_{(i_r,t)}\), and its
  initial term is a product of variables correspond to arcs from
  \(V^{+}\) to \(V^{-}\).
\end{description}
\begin{figure}[h]
  \begin{center}
    \epsfig{file=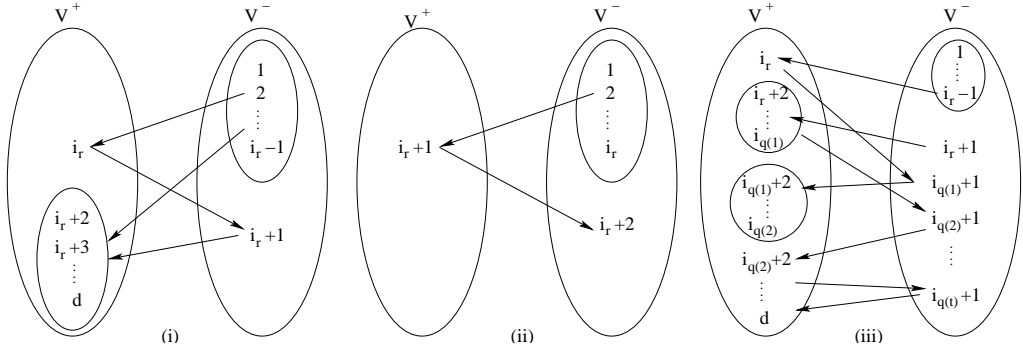,width=13.5cm}
    \caption{Cutsets corresponding to binomials (i), (ii) and (iii).}
  \end{center}
\end{figure}

If the case (i) occurs, the initial term of \(g^{(1)}_{(i_r)}\) is
\(x_{i_r,i_r+1}\), which contradicts that \(i_r\notin M(0)\). On the
other hand, a binomial of type (ii) can be reduced by \(g^{(1)}_{(1)}\)
by the above lemma, and cannot be contained in \({\mathcal G}\).

Let us consider that the case (iii) occurs. If
\(\vx^{\overline{T^{2}}}\) can be reduced by \(g^{(2)}_{(i_r,t)}\) with
\(t=1\). Then the monomial to which \(\vx^{\overline{T^{2}}}\) are
reduced by \(g^{(2)}_{(i_r,1)}\) can be reduced by a binomial \(f_{D}\in
I_{(I\ -M\T)}\), for the cutset \(D\) which corresponds to
\((V_{D}^{+},V_{D}^{-})\) where \(V_{D}^{+}=\{1,2,\ldots,i_r-1,i_r+1\}\)
and \(V_{D}^{-}=V\setminus V_{D}^{+}\), to some monomial \(m\) (the
initial term of \(f_{D}\) is a product of variables corresponds to arcs
from \(V_{D}^{+}\) to \(V_{D}^{-}\) since \(i_r\notin M(0)\)). 

 \begin{table}[h]
   \begin{center}
     \begin{tabular}{|l|l|l|l|}
       \hline
       \multicolumn{2}{|c|}{\em reduce\/ by\/ \(g^{(2)}_{(i_r,1)}\)\/} & 
       \multicolumn{2}{c|}{\em reduce\/ by\/ \(f_{D}\)\/}\\
       \hline
       divided variables & multiplied variables & divided variables &
       multiplied variable\\
       \hline
       \(x_{i_r,i_r+1}, x_{i_r,i_{q(1)}+1},\) & \(\{x_{k,i_r}:\ k\leq
       i_r-1\},\) &  \(\{x_{k,i_r}:  k\leq i_r-1\},\) &
       \(x_{i_r,i_r+1}\)\\
       \(x_{i_r+2,i_{q(1)}+1},\) & \(\{x_{k,l}:  k\leq i_r+1,\ k\neq
       i_r,\) &  \(\{x_{k,l}: k\leq i_r+1,\) & \\
       \(x_{i_r+3,i_{q(1)}+1},\) & \(\quad l\geq i_r+2, l\neq
       i_{q(1)}+1\},\) & \(\quad k\neq i_r, l\geq i_r+2\}\) & \\
       \(\ldots, x_{i_{q(1)},i_{q(1)}+1}\) & \(\{x_{i_{q(1)}+1,l}:\
       l\geq i_{q(1)}+2\}\) & & \\
       \hline
     \end{tabular}
   \end{center}
   \label{table:del_in_vars2}
   \caption{Divided and multiplied variables while reducing by
     \(g^{(2)}_{(i_r,1)}\) and \(f_{D}\).}
 \end{table}
 
 For a binomial \(f_{D'}\in I_{(I\ -M\T)}\) which corresponds to the
 cutset \(D'\) for \((V_{D'}^{+},V_{D'}^{-})\) such that
 \(V_{D'}^{-}=\{i_{q(1)}+1\}\) and \(V_{D'}^{+}=V\setminus V_{D'}^{-}\), 
 \(in_{\widetilde{\vb}}(f_{D'})\) corresponds to arcs from
 \(V_{D'}^{-}\) to \(V_{D'}^{+}\) (otherwise,
 \(\vx^{\overline{T_{\sigma\setminus\{i_r\}}}}\) can be reduced by
 \(f_{D'}\), which contradicts the assumption of the induction).  Then
 \(m\) can be reduced by \(f_{D'}\), and the resulting monomial is
 \(\vx^{\overline{T^{2}}}\)~(see Table 4), which contradicts to the
 definition of a term order by \(\widetilde{\vb}\).

Similarly, in the case that \(\vx^{\overline{T^{2}}}\) can be reduced by 
\(g^{(2)}_{(i_r,t)}\) for some \(t>1\), using the same \(f_{D}\) and
\(f_{D'}\in I_{(I\ -M\T)}\) which corresponds to the cutset \(D'\) for
\((V_{D'}^{+},V_{D'}^{-})\) such that
\(V_{D'}^{-}=\{i_{q(1)}+1,i_{q(2)}+1,\ldots,i_{q(t)}+1\}\), and
\(V_{D'}^{+}=V\setminus V_{D'}^{-}\), we can show the contradiction.

Therefore, \(\vx^{\overline{T^{2}}}\notin
in_{\widetilde{\vb}}(I_{(I\ -M\T)})\), and \(T^{2}\) is a desired
spanning tree \(T_{\sigma}\).
\end{proof}

\section*{Acknowledgment}
The authors thank Fumihiko Takeuchi, Masahiro Hachimori, Ryuichi
Hirabayashi, Kenji Kashiwabara, Masafumi Itoh, Takayuki Hibi and
Hidefumi Ohsugi for their useful comments.

\end{document}